\documentclass{amsart}

\newtheorem{theorem}{Theorem}[section]
\newtheorem{lemma}[theorem]{Lemma}
\newtheorem{corollary}[theorem]{Corollary}
\newtheorem{proposition}[theorem]{Proposition}

\theoremstyle{definition}
\newtheorem{definition}[theorem]{Definition}

\theoremstyle{remark}

\numberwithin{equation}{section}

\newcommand{\cB}{\mathcal{B}}
\newcommand{\cK}{\mathcal{K}}
\newcommand{\cM}{\mathcal{M}}
\newcommand{\cP}{\mathcal{P}}
\newcommand{\cR}{\mathcal{R}}
\newcommand{\C}{\mathbb{C}}
\newcommand{\N}{\mathbb{N}}
\newcommand{\R}{\mathbb{R}}
\newcommand{\T}{\mathbb{T}}
\newcommand{\Z}{\mathbb{Z}}
\newcommand{\eps}{\varepsilon}
\begin{document}

\title[Singular integral operators]{Fredholmness of Singular Integral Operators\\
with Piecewise Continuous Coefficients\\
on Weighted Banach Function Spaces}

\author{Alexei Yu. Karlovich}
\address{Departamento de Mathem\'atica,
Instituto Superior T\'ecnico,
Av. Rovisco Pais,
1049-001 Lisboa,
Portugal}
\email{akarlov@math.ist.utl.pt}

\thanks{The author is partially supported by F.C.T. (Portugal)
grants POCTI 34222/MAT/2000 and PRAXIS XXI/BPD/22006/99.}

%    General info

\date{February 18, 2003}
%%%%%%%%%%%%%%%%%%%%%%%%%%%%%%%%%%%%%%%%%%%%%%%%%%%%%%%%%%%%%%%%%%%
\begin{abstract}
We prove necessary conditions for the Fredholmness of singular integral
operators with piecewise continuous coefficients on weighted Banach
function spaces. These conditions are formulated in terms of
indices of submultiplicative functions associated with local
properties of the space, of the curve, and of the weight.
As an example, we consider weighted Nakano spaces
$L^{p(\cdot)}_w$  (weighted Lebesgue spaces with variable
exponent).  Moreover, our necessary conditions become also sufficient
for weighted Nakano spaces over nice curves whenever
$w$ is a Khvedelidze weight, and the variable exponent $p(t)$
satisfies the estimate $|p(\tau)-p(t)|\le A/(-\log|\tau-t|)$.
\end{abstract}
\maketitle
%%%%%%%%%%%%%%%%%%%%%%%%%%%%%%%%%%%%%%%%%%%%%%%%%%%%%%%%%%%%%%%%%%%
\section{Introduction}
Let $\Gamma$ be a Jordan curve, that is, a curve that  
homeomorphic to a circle. We suppose that $\Gamma$ is rectifiable. 
We equip $\Gamma$ with Lebesgue length measure $|d\tau|$ and
the counter-clockwise orientation. The \textit{Cauchy singular integral}
of a measurable function $f:\Gamma\to\C$ is defined by
\[
(Sf)(t):=\lim_{R\to 0}\frac{1}{\pi i}\int_{\Gamma\setminus\Gamma(t,R)}
\frac{f(\tau)}{\tau-t}d\tau
\quad (t\in\Gamma),
\]
where the ``portion'' $\Gamma(t,R)$ is
\[
\Gamma(t,R):=\{\tau\in\Gamma:|\tau-t|<R\}
\quad (R>0).
\]
It is well known that $(Sf)(t)$ exists a.e. on $\Gamma$
whenever $f$ is integrable (see \cite[Theorem~2.22]{Dynkin87}).
A measurable function $w:\Gamma\to[0,\infty]$ is referred to as a
\textit{weight} if $0<w(t)<\infty$ a.e. on $\Gamma$.
The Cauchy singular integral generates a bounded linear operator
$S$ on the weighted Lebesgue space $L^p_w (1<p<\infty)$ with the norm
\[
\|f\|_{L^p_w}:=\left(\int_\Gamma |f(\tau)|^pw^p(\tau)|d\tau|\right)^{1/p}
\]
if and only if $w$ is a Muckenhoupt weight ($w\in A_p(\Gamma)$), that is,
\[
\sup_{t\in\Gamma}\sup_{R>0}
\left(\frac{1}{R}\int_{\Gamma(t,R)}w^p(\tau)|d\tau|\right)^{1/p}
\left(\frac{1}{R}\int_{\Gamma(t,R)}w^{-p'}(\tau)|d\tau|\right)^{1/p'}<\infty,
\frac{1}{p}+\frac{1}{p'}=1
\]
(see, e.g., \cite[Theorem~4.15]{BK97}).
By H\"older's inequality, if $w\in A_p(\Gamma)$, then $\Gamma$ is a
\textit{Carleson} (or \textit{Ahlfors-David regular}) \textit{curve},
that is,
%%%
\begin{equation}\label{eq:Carleson}
C_\Gamma:=\sup_{t\in\Gamma}\sup_{R>0}\frac{|\Gamma(t,R)|}{R}<\infty,
\end{equation}
%%%
where $|\Omega|$ denotes the measure of a measurable set $\Omega\subset\Gamma$.
The constant $C_\Gamma$ is said to be the \textit{Carleson constant}.
We denote by $PC$ the Banach algebra of all \textit{piecewise continuous}
functions on the curve $\Gamma$: by definition, $a$ is in $PC$
if and only if $a$ is in $L^\infty$ and the one-sided limits
\[
a(t\pm 0):=\lim_{\tau\to t\pm 0} a(\tau)
\]
exist for every $t\in\Gamma$.

A bounded linear operator $A$ on a Banach space is said to be
\textit{semi-Fredholm} if its image is closed and at
least one of the so-called defect numbers
\[
n(A):=\dim\ker A, \quad d(A):=\dim\ker A^*
\]
is finite. A semi-Fredholm operator $A$ is called \textit{Fredholm} if both
$n(A)$ and $d(A)$ are finite. In this case the difference $n(A)-d(A)$ is
referred to as the \textit{index}  of the operator $A$. Basic properties
of (semi)-Fredholm operators are  discussed in \cite{BS90,GK92,MP86} and in
many other monographs.

The study of Fredholmness of one-dimensional singular integral operators
of the form
\[
R_a:=aP_++P_-,\quad
a\in PC,\quad
P_\pm:=(I\pm S)/2
\]
on Lebesgue spaces with power (Khvedelidze) weights
%%%
\begin{equation}\label{eq:Khvedelidze-def}
\varrho(t):=\prod_{k=1}^n |t-\tau_k|^{\lambda_k},
\quad \tau_k\in\Gamma,\quad k\in\{1,\dots, n\},\quad n\in\N,
\end{equation}
%%%
over Lyapunov curves started in the fiftieth with B.~V.~Khvedelidze
\cite{Khvedelidze56} and was continued in the sixties by
H.~Widom, I.~B.~Simonenko, I.~Gohberg and N.~Krupnik, and others.
The history and corresponding references can be found, e.g., in
\cite{BK97,GK92,KS01,Khvedelidze75,MP86}.
In the beginning of nineties, I.~Spitkovsky proved Fredholm criteria
for singular integral operators with piecewise continuous coefficients
on Lebesgue spaces with Muckenhoupt weights over smooth curves
\cite{Spitkovsky92}. In the middle of nineties, A.~B\"ottcher
and Yu.~I.~Karlovich accomplished the Fredholm theory for the
algebra of singular integral operators with piecewise continuous
coefficients on Lebesgue spaces with Muckenhoupt weights over
general Carleson curves. These results are documented in
\cite{BK97}; see also the brief but nice presentation in \cite{BK01}.

Lebesgue spaces $L^p,1\le p\le\infty$, are the simplest examples
of so-called \textit{Banach function spaces} introduced by
W.~A.~J.~Luxemburg in 1955. This scale of spaces includes
Orlicz, Lorentz, and all other \textit{rearrangement-invariant spaces}.
By analogy with weighted Lebesgue spaces, for a Banach function space
$X$ and a  weight $w$, it is possible to define the
\textit{weighted Banach function space}
\[
X_w:=\Big\{f\ \mbox{is measurable on }\Gamma \ \mbox{and} \ fw\in X\Big\}.
\]
Under some restrictions on the weight $w$, the space $X_w$ is itself
a Banach function space, although if $X$ is a rearrangement-invariant
Banach function space, then $X_w$ is not necessarily
rearrangement-invariant (even if $X$ is a Lebesgue space). Another
interesting class of Banach function spaces which are not
rearrangement-invariant constituted by Nakano spaces $L^{p(\cdot)}$
(generalized Lebesgue spaces with variable exponent).
For details and references, see Section~\ref{sec:WBFS}.

Unfortunately, few is known about the boundedness of $S$
on general weighted Banach function spaces $X_w$. As far as we know,
even a criterion for the boundedness of $S$ on Orlicz spaces
$L^\varphi_w$ with general weights $w$ over general Carleson curves
is unknown at the moment (February of 2003). We proved necessary
conditions for the boundedness of $S$ on weighted
rearrangement-invariant Banach function spaces \cite[Theorem~3.2]{K98}
in terms of an analog of the Muckenhoupt class. On the other hand,
if a weight $w$ belongs to the Muckenhoupt classes
$A_{1/\alpha_X}(\Gamma)$ and $A_{1/\beta_X}(\Gamma)$
where $\alpha_X,\beta_X\in (0,1)$ are the Boyd indices of a
rearrangement-invariant Banach function space $X$, then $S$ is
bounded on the weighted rearrangement-invariant Banach function
space $X_w$ (see \cite[Theorem~4.5]{K02}).

On the basis of these boundedness results, following the approach of
A.~B\"ott\-cher, Yu.~Karlovich, and I.~Spitkovsky, the author
proved separately necessary and sufficient conditions
for Fredholmness of singular integral operators with piecewise
continuous coefficients on weighted rearrangement-invariant Banach
function spaces \cite{K00,K02}. Under some restrictions on
spaces, curves, and weights, these conditions coincide, that is,
become criteria. In that cases, the Banach algebra of singular
integral operators with piecewise continuous coefficients is also
studied \cite{K02}.

Very recently V.~M.~Kokilashvili and S.~G.~Samko have proved
criteria for the boundedness of $S$ on Nakano spaces
$L^{p(\cdot)}_\varrho$ with Khvedelidze weights $\varrho$
over Lyapunov curves or Radon curves without cusps provided the variable
exponent $p$ satisfies the estimate
%%%
\begin{equation}\label{eq:p-restriction}
|p(\tau)-p(t)|\le A/(-\log|\tau-t|),
\quad \tau,t\in\Gamma,
\quad |\tau-t|\le 1/2
\end{equation}
%%%
(see \cite[Theorem~2]{KS02-3} or Theorem~\ref{th:Samko-Kokilashvili}).
With the help of this key result,
they have proved Fredholm criteria for the operator $aP_++bP_-$
with piecewise continuous functions $a,b$ having finite numbers
of jumps on (non-weighted) Nakano spaces $L^{p(\cdot)}$
(see \cite[Theorem~A]{KS03-4}).

For an arbitrary weight $w$ and an arbitrary Banach function
space $X$, we define the weighted Banach function space $X_w$.
Assume that
%%%
\begin{enumerate}
\item[{\rm (B)}] the Cauchy singular integral operator $S$
is bounded on $X_w$;
\item[{\rm (R)}] $X_w$ is reflexive.
\end{enumerate}
%%%
We show that property (B) implies the condition $A_X(\Gamma)$
of Muckenhoupt type. In that case $X_w$ is itself
a Banach function space. Under the assumptions (B) and (R)
we prove necessary conditions for Fredholmness of
singular integral operators $R_a$ with piecewise continuous
coefficients $a$ in the weighted Banach function spaces $X_w$.
This result generalizes corresponding
necessary conditions in \cite[Theorem~4.2]{K00}. As an example, we
consider these necessary conditions in Nakano spaces
$L^{p(\cdot)}_w$ with general weights $w$. They have almost
the same form as in the case of Lebesgue spaces $L^p_w$ with Muckenhoupt
weights over Carleson curves (see \cite[Proposition~7.3]{BK97}).
We should only replace the constant $p$ (for weighted Lebesgue
spaces $L^p_w$) by the value $p(t)$ of the variable exponent $p(\cdot)$
at each point $t\in\Gamma$ (for weighted Nakano spaces $L^{p(\cdot)}_w$).
Our approach is based on a local principle of Simonenko type, the
Wiener-Hopf factorization of local representatives, and the
theory of submultiplicative functions associated with local
properties of the curve, of the weight, and of the space.
Using of the local principle allows us to consider coefficients
$a$ having a countable number of jumps (in contrast to
\cite{KS03-4}, where only a finite number of jumps is allowed).

The paper is organized as follows. In Section~\ref{sec:WBFS}
we collect necessary preliminaries on weighted Banach function
spaces $X_w$ and Nakano spaces $L^{p(\cdot)}$. In Section~\ref{sec:AX}
we define an analog of the Muckenhoupt class $A_p(\Gamma)$,
replacing the norm in $L^p$ by the norm in a Banach function
space $X$. We denote this class by $A_X(\Gamma)$. We show that
if $w\in A_X(\Gamma)$ and $1\in A_X(\Gamma)$, then
$\log w$ has bounded mean oscillation. In Section~\ref{sec:indices}
we remind the definitions and some properties of submultiplicative functions
associated with the local behavior of the curve, of the weight, and
of the space. In Section~\ref{sec:relations} we study inequalities
between the indices of submultiplicative functions defined
in Section~\ref{sec:indices}. We investigate so-called
\textit{indicator functions} $\alpha_t^*,\beta_t^*$
and $\alpha_t,\beta_t$ of the triple $(\Gamma,X,w)$ and of the
pair $(\Gamma,w)$, respectively. In particular, we show that if
$X$ is a Nakano space $L^{p(\cdot)}$ with a variable exponent
$p(\cdot)$ satisfying (\ref{eq:p-restriction}),
then we can separate the influence of the space from the influence
of the weight and the curve, that is,
$\alpha_t^*(x)=1/p(t)+\alpha_t(x), \beta_t^*(x)=1/p(t)+\beta_t(x)$
for $x\in\R$ such that
$|(\tau-t)^{y+ix}|w(\tau)\in A_{L^{p(\cdot)}}(\Gamma,t)$,
where $A_{L^{p(\cdot)}}(\Gamma,t)$ is the local analog of
$A_{L^{p(\cdot)}}(\Gamma)$.
So, weighted Nakano spaces
satisfy the ``disintegration condition'' in the terminology of
\cite{K98,K02}.

In Section~\ref{sec:bounded} we prove that the condition
$w\in A_X(\Gamma)$ is necessary
for the boundedness of the Cauchy singular integral operator
$S$ on the weighted Banach function spaces $X_w$.
Further we extend basic results on the Fredholmness
of singular integral operators with bounded measurable coefficients
(the local principle, the theorem about a Wiener-Hopf factorization, etc.)
to weighted Banach function spaces satisfying Axioms (B) and (R).
These results are natural extensions of the classical theory for
Lebesgue spaces with Khvedelidze weights over Lyapunov curves
(see, e.g., \cite[Ch.~7-8]{GK92} or \cite[Ch.~4]{MP86}).
A canonical local representative $g_{t,\gamma} \ (t\in\Gamma,\gamma\in\C)$
for a piecewise continuous function is constructed in Section~\ref{sec:PC}.
We prove separately necessary and sufficient conditions for factorability of
$g_{t,\gamma}$ in the weighted Banach function space $X_w$. On the basis of our
necessary conditions for factorability, with the help of the results
of Section~\ref{sec:bounded}, we prove necessary conditions for
Fredholmness of the singular integral operator $R_a=aP_++P_-$
with $a\in PC$ in $X_w$. These conditions are formulated in terms
of the indicator functions $\alpha_t^*$ and $\beta_t^*$ defined
in Section~\ref{sec:relations}. In Section~\ref{sec:PC-Nakano}
we reformulate these necessary conditions for weighted Nakano spaces
$L^{p(\cdot)}_w$ with general weights $w$ and variable exponents
satisfying (\ref{eq:p-restriction}) in terms of simpler indicator
functions $\alpha_t$ and $\beta_t$. With the help of the boundedness
criteria by V.~M.~Kokilashvili and S.~G.~Samko \cite[Theorem~2]{KS02-3},
we prove that the latter necessary conditions become also
sufficient if  $w=\varrho$ is a Khvedelidze weight and $\Gamma$ is either
a Lyapunov Jordan curve or a Radon Jordan curve without cusps.
%%%%%%%%%%%%%%%%%%%%%%%%%%%%%%%%%%%%%%%%%%%%%%%%%%%%%%%%%%%%%%%%%%%%%%%%
\section{Weighted Banach function spaces}\label{sec:WBFS}
\subsection{Banach function spaces}
Let $\Gamma$ be a  rectifiable Jordan (i.e., homeomorphic to a circle)
curve equipped with Lebesgue length measure $|d\tau|$.
The set of all measurable complex-valued functions on $\Gamma$ is
denoted by $\cM$. Let $\cM^+$ be the subset of functions in $\cM$ whose
values lie  in $[0,\infty]$. The characteristic function of a measurable
set $E\subset\Gamma$ is denoted by $\chi_E$.
%%%%%%%%%%%%%%%%%%%%%%%%%%%%%%%%%%%%%%%%%%%%%%%%%%%%%%%%%%%%%%%%%%%%%%%%
\begin{definition}\label{def:Luxemburg}
(\textbf{W.~A.~J.~Luxemburg, 1955}, see \cite[Ch.~1, Definition~1.1]{BS88}).
A mapping $\rho:\cM^+\to [0,\infty]$ is called a {\it Banach function norm}
if, for all functions $f,g, f_n \ (n\in\N)$ in $\cM^+$, for all constants
$a\ge 0$, and for all measurable subsets $E$ of $\Gamma$, the
following properties hold:
%%%%
%\newpage
\begin{eqnarray*}
{\rm (A1)} & &
\rho(f)=0  \Leftrightarrow  f=0\ \mbox{a.e.},
\quad
\rho(af)=a\rho(f),
\quad
\rho(f+g) \le \rho(f)+\rho(g),\\
{\rm (A2)} & &0\le g \le f \ \mbox{a.e.} \ \Rightarrow \ \rho(g) \le \rho(f)
\quad\mbox{(the lattice property)},\\
{\rm (A3)} & &0\le f_n \uparrow f \ \mbox{a.e.} \ \Rightarrow \
       \rho(f_n) \uparrow \rho(f)\quad\mbox{(the Fatou property)},\\
{\rm (A4)} & &\rho(\chi_E) <\infty,\\
{\rm (A5)} & &\int_E f(\tau)|d\tau| \le C_E\rho(f)
\end{eqnarray*}
%%%%
with $C_E \in (0,\infty)$ may depend on $E$ and $\rho$ but is independent of $f$.
\end{definition}
%%%%%%%%%%%%%%%%%%%%%%%%%%%%%%%%%%%%%%%%%%%%%%%%%%%%%%%%%%%%%%%%%%%%%%%%
When functions differing only on a set of measure zero are identified,
the set $X$ of all functions $f\in\cM$ for which $\rho(|f|)<\infty$ is
called a \textit{Banach function space}. For each $f\in X$, the norm of
$f$ is defined by
\[
\|f\|_X :=\rho(|f|).
\]
The set $X$ under the natural linear space operations and under this norm
becomes a Banach space (see \cite[Ch.~1, Theorems~1.4 and~1.6]{BS88}).

If $\rho$ is a Banach function norm, its associate norm $\rho'$ is
defined on $\cM^+$ by
\[
\rho'(g):=\sup\left\{
\int_\Gamma f(\tau)g(\tau)|d\tau| \ : \ f\in \cM^+, \ \rho(f) \le 1
\right\}, \quad g\in \cM^+.
\]
It is a Banach function norm itself \cite[Ch.~1, Theorem~2.2]{BS88}.
The Banach function space $X'$ determined by the Banach function norm
$\rho'$ is called the {\it associate space (K\"othe dual)} of $X$.
The associate space $X'$ is a subspace of the dual space $X^*$.
The construction of the associate space implies the following
H\"older inequality for Banach function spaces.
%%%%%%%%%%%%%%%%%%%%%%%%%%%%%%%%%%%%%%%%%%%%%%%%%%%%%%%%%%%%%%%%%%%%%%%%%%%
\begin{lemma}\label{le:Hoelder}
{\rm (see \cite[Ch.~1, Theorem~2.4]{BS88}).}
Let $X$ be a Banach function  space and $X'$ be its associate space.
If $f\in X$ and $g\in X'$, then $fg$ is integrable and
\[
\|fg\|_{L^1}\le \|f\|_X\|g\|_{X'}.
\]
\end{lemma}
%%%%%%%%%%%%%%%%%%%%%%%%%%%%%%%%%%%%%%%%%%%%%%%%%%%%%%%%%%%%%%%%%%%%%%%%%%%
\subsection{Rearrangement-invariant Banach function spaces}
Let $\cM_0$ and $\cM_0^+$ be the classes of a.e. finite
functions in $\cM$ and $\cM^+$, respectively. Two functions
$f,g\in\cM_0$ are said to be equimeasurable if
\[
\Big|\{\tau\in\Gamma:|f(\tau)|>\lambda\}\Big|=
\Big|\{\tau\in\Gamma:|g(\tau)|>\lambda\}\Big|
\quad\mbox{for all}\quad \lambda\ge 0.
\]
A Banach function norm $\rho:\cM^+ \to [0,\infty]$ is called
rearrangement-invariant if for every pair of equimeasurable
functions $f,g \in \cM^+_0$ the equality  $\rho(f)=\rho(g)$
holds. In that case, the Banach function space $X$ generated
by $\rho$ is said to be a {\it rearrangement-invariant Banach function space}
(or simply rearrangement-invariant space).
Lebes\-gue, Orlicz, Lorentz, and Lorentz-Orlicz
spaces are classical examples of rearrange\-ment-invariant
Banach function spaces (see, e.g.,  \cite{BS88} and the references therein).

If $X$ is an arbitrary rearrangement-invariant Banach function space
and $X'$ is its associate space, then for a measurable set $E\subset\Gamma$,
%%%
\begin{equation}\label{eq:RI-fundamental}
\|\chi_E\|_X\|\chi_E\|_{X'}=|E|
\end{equation}
%%%
(see, e.g., \cite[Ch.~2, Theorem~5.2]{BS88}).
%%%%%%%%%%%%%%%%%%%%%%%%%%%%%%%%%%%%%%%%%%%%%%%%%%%%%%%%%%%%%%%%%%%%%%%%%%%
\subsection{Nakano spaces $L^{p(\cdot)}$}
Function spaces $L^{p(\cdot)}$ of Lebesgue type with variable exponent
$p$ were studied for the first time probably  by W.~Orlicz  \cite{Orlicz31} in 1931.
Inspired by the successful theory of Orlicz spaces, H.~Nakano defined in the late forties
\cite{Nakano50,Nakano51} so-called \textit{modular spaces}. He considered the space
$L^{p(\cdot)}$ as an example of modular spaces. J.~Musielak and W.~Orlicz \cite{MO59}
in 1959 extended the definition of modular spaces by H.~Nakano.
Actually, that paper was the starting point for the theory of Musielak-Orlicz
spaces (generalized Orlicz spaces generated by Young
functions with a parameter), see \cite{Musielak83}.

Let $p:\Gamma\to[1,\infty)$ be a measurable function. Consider the
convex modular (see \cite[Ch.~1]{Musielak83} for definitions and properties)
\[
m(f,p):=\int_\Gamma|f(\tau)|^{p(\tau)}|d\tau|.
\]
Denote by $L^{p(\cdot)}$ the set of all measurable complex-valued functions
$f$ on $\Gamma$ such that $m(\lambda f,p)<\infty$ for some $\lambda=\lambda(f)>0$.
This set becomes a Banach space with respect to the \textit{Luxemburg-Nakano norm}
\[
\|f\|_{L^{p(\cdot)}}:=\inf\Big\{\lambda>0: \ m(f/\lambda,p)\le 1\Big\}
\]
(see, e.g., \cite[Ch.~2]{Musielak83}). So, the spaces $L^{p(\cdot)}$
are a special case of Musielak-Orlicz spaces. Sometimes the spaces
$L^{p(\cdot)}$ are referred to as Nakano spaces (see, e.g.,
\cite[p.~151]{FJK92}, \cite[p.~179]{KT90}).
We will follow this tradition. Clearly, if $p(\cdot)=p$ is constant,
then the Nakano space $L^{p(\cdot)}$ is isometrically isomorphic
to the Lebesgue space $L^p$. Therefore, sometimes $L^{p(\cdot)}$
are called generalized Lebesgue spaces with variable exponent.
%%%%%%%%%%%%%%%%%%%%%%%%%%%%%%%%%%%%%%%%%%%%%%%%%%%%%%%%%%%%%%%%%%%%%%%%%%%
\begin{lemma}\label{le:Nakano-BFS}
{\rm (see, e.g., \cite[Proposition~1.3]{ELN99})}.
Let $p:\Gamma\to[1,\infty)$ be a measurable function. The Nakano space
$L^{p(\cdot)}$ is a Banach function space.
\end{lemma}
%%%%%%%%%%%%%%%%%%%%%%%%%%%%%%%%%%%%%%%%%%%%%%%%%%%%%%%%%%%%%%%%%%%%%%%%%%%
It is not difficult to show that $L^{p(\cdot)}$ is not rearrangement-invariant,
in general.

The following result on the reflexivity and duality of Nakano spaces
was precisely stated in \cite[Theorem~2.3 and Corollary~2.7]{KR91},
although it can be obtained from more general results for Musielak-Orlicz
spaces \cite[Ch.~1--2]{Musielak83} (see also \cite{Orlicz31}).
%%%%%%%%%%%%%%%%%%%%%%%%%%%%%%%%%%%%%%%%%%%%%%%%%%%%%%%%%%%%%%%%%%%%%%%%%%%
\begin{lemma}\label{le:Nakano-duality}
Let $p:\Gamma\to[1,\infty)$ be a measurable function.
If
\[
1<\operatornamewithlimits{ess\,inf}_{t\in\Gamma} p(t)
\le
\operatornamewithlimits{ess\,sup}_{t\in\Gamma} p(t)<\infty,
\]
then the Nakano space $L^{p(\cdot)}$ is reflexive. Its associate space
coincides (up to the equivalence of the norms) with the Nakano space
$L^{p'(\cdot)}$, where
\[
p'(\tau):=\frac{p(\tau)}{p(\tau)-1}.
\]
\end{lemma}
%%%%%%%%%%%%%%%%%%%%%%%%%%%%%%%%%%%%%%%%%%%%%%%%%%%%%%%%%%%%%%%%%%%%%%%%%%%
Finally, Nakano spaces are important in applications to fluid dynamics
\cite{Ruzicka00}.
%%%%%%%%%%%%%%%%%%%%%%%%%%%%%%%%%%%%%%%%%%%%%%%%%%%%%%%%%%%%%%%%%%%%%%%%%%%
\subsection{Weighted Banach function spaces}\label{subsection:WBFS}
Let $X$ be a Banach function space generated by a Banach function norm
$\rho$ and let $w:\Gamma\to[0,\infty]$ be a weight. Define the mapping
$\rho_w:\cM^+\to[0,\infty]$ and the set $X_w$ by
\[
\rho_w(f):=\rho(fw)
\quad (f\in\cM^+),
\quad\quad
X_w:=\Big\{f\in\cM^+:\quad fw\in X\Big\}.
\]
%%%%%%%%%%%%%%%%%%%%%%%%%%%%%%%%%%%%%%%%%%%%%%%%%%%%%%%%%%%%%%%%%%%%%%%%%%
\begin{lemma}\label{le:WBFS}
{\rm (a)} $\rho_w$ satisfies Axioms {\rm (A1)--(A3)} in
Definition~{\rm\ref{def:Luxemburg}} and $X_w$ is a linear normed space with
respect to the norm
\[
\|f\|_{X_w}:=\rho_w(|f|)=\rho(|fw|)=\|fw\|_X;
\]

\noindent
{\rm (b)} if $w\in X$ and $1/w\in X'$, then $\rho_w$ is a Banach function
norm and $X_w$ is a Banach function space generated by $\rho_w$.
Moreover,
\[
L^\infty\subset X_w\subset L^1;
\]

\noindent
{\rm (c)} if $w\in X$ and $1/w\in X'$, then $X_{1/w}'$ is the associate
space for the Banach function space $X_w$.
\end{lemma}
%%%%%%%%%%%%%%%%%%%%%%%%%%%%%%%%%%%%%%%%%%%%%%%%%%%%%%%%%%%%%%%%%%%%%%%%%
\begin{proof}
Part (a) follows from Axioms (A1)--(A3) for the Banach function norm $\rho$
and the fact that $0<w(\tau)<\infty$ a.e. on $\Gamma$.

(b) If $w\in X$, then by Axiom (A2) for $\rho$, we get $w\chi_E\in X$
for every measurable set $E$ of $\Gamma$. Therefore,
$\rho_w(\chi_E)=\rho(w\chi_E)<\infty$. Thus, $\rho_w$ satisfies
Axiom (A4). By H\"older's inequality (see Lemma~\ref{le:Hoelder})
and Axiom (A2) for $\rho$, we have
%%%
\begin{eqnarray}
\label{eq:WBFS-1}
\int_E f(\tau)|d\tau|
&=&
\int_\Gamma \Big(f(\tau)w(\tau)\chi_E(\tau)\Big)\frac{\chi_E(\tau)}{w(\tau)}|d\tau|
\\
&\le&
\rho(fw\chi_E)\rho'(\chi_E/w)\le\rho(fw)\rho'(\chi_E/w)=:C_E\rho_w(f),
\nonumber
\end{eqnarray}
%%%
where $C_E:=\rho'(\chi_E/w)\in(0,\infty)$. This constant, clearly,
depends on $\rho,w$ (and thus on $\rho_w$) and $E$, but it is independent
of $f$. Therefore, $\rho_w$ satisfies Axiom (A5). Thus, $\rho_w$
is a Banach function norm and $X_w$ is a Banach function space.

From (\ref{eq:WBFS-1}) and Axiom (A2) for $X'$ it follows that
\[
\|f\|_{L^1}\le \|f\|_{X_w}\|1/w\|_{X'},\quad f\in X_w.
\]
Hence, $X_w\subset L^1$, in view of $1/w\in X'$. On the other hand,
for $f\in L^\infty$,
\[
0\le |f(\tau)|\le \|f\|_\infty
\quad\mbox{a.e. on}\quad\Gamma.
\]
By Axioms
(A2) and (A1) for $\rho_w$, we have
\[
\|f\|_{X_w}=\rho_w(|f|)\le\rho_w(\|f\|_\infty)=\|f\|_\infty\rho_w(1)
=\|f\|_\infty\|w\|_X.
\]
Thus, $L^\infty\subset X_w$, in view of $w\in X$. Part (b) is proved.

(c) For $g\in\cM^+$, we have
%%%
\begin{eqnarray*}
(\rho_w)'(g)
&=&
\sup\left\{\int_\Gamma f(\tau)g(\tau)|d\tau|:
\quad f\in\cM^+,\quad \rho_w(f)\le 1
\right\}
\\
&=&
\sup\left\{\int_\Gamma \Big(f(\tau)w(\tau)\Big)
\left(\frac{g(\tau)}{w(\tau)}\right)|d\tau|:
\quad f\in\cM^+,\quad \rho(fw)\le 1
\right\}
\\
&=&
\sup\left\{\int_\Gamma h(\tau)
\left(\frac{g(\tau)}{w(\tau)}\right)|d\tau|:
\quad h\in\cM^+,\quad \rho(h)\le 1
\right\}
\\
&=&\rho'(g/w).
\end{eqnarray*}
%%%
Hence, $(X_w)'=X_{1/w}'$.
\end{proof}
%%%%%%%%%%%%%%%%%%%%%%%%%%%%%%%%%%%%%%%%%%%%%%%%%%%%%%%%%%%%%%%%%%%%%%%%%

We will refer to the normed space $X_w$ as a \textit{weighted Banach
function space} generated by the Banach function space $X$ and the weight
$w$. From Lemma~\ref{le:WBFS}(b) it follows that the weighted Banach
function space $X_w$ is a Banach function space itself whenever $w\in X$
and $1/w\in X'$.

For other definition (different from our) of weighted Banach function
spaces, see, e.g., \cite{KOPR93,LNP03}.
%%%%%%%%%%%%%%%%%%%%%%%%%%%%%%%%%%%%%%%%%%%%%%%%%%%%%%%%%%%%%%%%%%%%%%%%%%%
\subsection{Separability and reflexivity of weighted Banach function spaces}
%%%
A function $f$ in a Banach function space $X$ is said to have
\textit{absolutely continuous norm} in $X$ if $\|f\chi_{E_n}\|_X\to 0$
for every sequence $\{E_n\}_{n=1}^\infty$ of measurable sets on $\Gamma$
satisfying $\chi_{E_n}\to 0$ a.e. on $\Gamma$ as $n\to\infty$. If all
functions $f\in X$ have this property, then the space $X$ itself is said to have
\textit{absolutely continuous norm} (see \cite[Ch.~1, Section~3]{BS88}).

In this subsection we assume that $X$ is a Banach function space and
$w$ is a weight such that $w\in X$ and $1/w\in X'$. Then, by
Lemma~\ref{le:WBFS}(b), the weighted Banach function space $X_w$ is itself
a Banach function space.
%%%%%%%%%%%%%%%%%%%%%%%%%%%%%%%%%%%%%%%%%%%%%%%%%%%%%%%%%%%%%%%%%%%%%%%%%%%
\begin{proposition}\label{pr:AC}
If $X$ has absolutely continuous norm, then $X_w$ has absolutely
continuous norm too.
\end{proposition}
%%%%%%%%%%%%%%%%%%%%%%%%%%%%%%%%%%%%%%%%%%%%%%%%%%%%%%%%%%%%%%%%%%%%%%%%%%%
\begin{proof}
If $f\in X_w$, then $fw\in X$ has absolutely continuous norm in $X$. Therefore,
$\|f\chi_{E_n}\|_{X_w}=\|fw\chi_{E_n}\|_X\to 0$ for every sequence
$\{E_n\}_{n=1}^\infty$ of measurable sets on $\Gamma$ satisfying
$\chi_{E_n}\to 0$ a.e. on $\Gamma$ as $n\to \infty$. Thus, $f\in X_w$
has absolutely continuous norm in $X_w$.
\end{proof}

From Lemma~\ref{le:WBFS} and \cite[Ch.~1, Corollaries 4.3, 4.4]{BS88}
we obtain the following.
%%%%%%%%%%%%%%%%%%%%%%%%%%%%%%%%%%%%%%%%%%%%%%%%%%%%%%%%%%%%%%%%%%%%%%%%%%%
\begin{lemma}\label{le:duality-reflexivity}
{\rm (a)}
The Banach space dual $(X_w)^*$ of the weighted Banach function space
$X_w$ is isometrically isomorphic to the associate space
$X_{1/w}'$ if and only if $X_w$ has absolutely continuous norm.
If $X_w$ has absolutely continuous norm, then the general form of a linear
functional on $X_w$ is given by
\[
G(f):=\int_\Gamma f(\tau)\overline{g(\tau)}|d\tau|,
\quad
g\in X_{1/w}',
\quad\mbox{and}\quad \|G\|_{(X_w)^*}=\|g\|_{X_{1/w}'}.
\]
{\rm (b)} The weighted Banach function space $X_w$ is reflexive if and only if
both $X_w$ and $X_{1/w}'$ have absolutely continuous norm.
\end{lemma}
%%%%%%%%%%%%%%%%%%%%%%%%%%%%%%%%%%%%%%%%%%%%%%%%%%%%%%%%%%%%%%%%%%%%%%%%%%%
\begin{corollary}\label{co:weighted-reflexivity}
If $X$ is reflexive, then $X_w$ is reflexive.
\end{corollary}
%%%%%%%%%%%%%%%%%%%%%%%%%%%%%%%%%%%%%%%%%%%%%%%%%%%%%%%%%%%%%%%%%%%%%%%%%%%
\begin{proof}
If $X$ is reflexive, then, by \cite[Ch.~1, Corollary~4.4]{BS88},
both $X$ and $X'$ have absolutely continuous norm. In that case, due to
Proposition~\ref{pr:AC}, both $X_w$ and $X_{1/w}'$ have absolutely
continuous norm. By Lemma~\ref{le:duality-reflexivity}(b), $X_w$ is reflexive.
\end{proof}
%%%%%%%%%%%%%%%%%%%%%%%%%%%%%%%%%%%%%%%%%%%%%%%%%%%%%%%%%%%%%%%%%%%%%%%%%%%
Since Lebesgue length measure $|d\tau|$ is separable (for the definition of
a separable measure, see, e.g., \cite[p.~27]{BS88} or
\cite[Section~6.10]{KA84}), from Lemma~\ref{le:WBFS} and
\cite[Ch.~1, Corollary~5.6]{BS88} we immediately get the following criterion.
%%%%%%%%%%%%%%%%%%%%%%%%%%%%%%%%%%%%%%%%%%%%%%%%%%%%%%%%%%%%%%%%%%%%%%%%%%%
\begin{lemma}\label{le:weighted-separability}
The weighted Banach function space $X_w$ is separable if and only if it has
absolutely continuous norm.
\end{lemma}
%%%%%%%%%%%%%%%%%%%%%%%%%%%%%%%%%%%%%%%%%%%%%%%%%%%%%%%%%%%%%%%%%%%%%%%%%%%
We denote by $C$ the set of all continuous functions on $\Gamma$
and by $\cR$ the set of all rational functions without poles on
the curve $\Gamma$. With the help of Lemmas~\ref{le:duality-reflexivity}
and~\ref{le:weighted-separability}, literally repeating  the proof of
\cite[Lemma~1.3]{K00}, one can get the following.
%%%%%%%%%%%%%%%%%%%%%%%%%%%%%%%%%%%%%%%%%%%%%%%%%%%%%%%%%%%%%%%%%%%%%%%%%%%
\begin{lemma}\label{le:C-density}
The weighted Banach function space $X_w$ is separable if and only if
$C$ is dense in $X_w$.
\end{lemma}
%%%%%%%%%%%%%%%%%%%%%%%%%%%%%%%%%%%%%%%%%%%%%%%%%%%%%%%%%%%%%%%%%%%%%%%%%%%
\begin{corollary}\label{co:R-density}
If $X_w$ (or $X$) is reflexive, then $\cR$ is dense in $X_w$ and in its
associate space $X_{1/w}'$.
\end{corollary}
%%%%%%%%%%%%%%%%%%%%%%%%%%%%%%%%%%%%%%%%%%%%%%%%%%%%%%%%%%%%%%%%%%%%%%%%%%%
\begin{proof}
If $X_w$ is reflexive, then by Lemmas~\ref{le:duality-reflexivity}(b)
and~\ref{le:weighted-separability}, both $X_w$ and $X_{1/w}'$
are separable. This implies that $C$ is dense in $X_w$ and in
$X_{1/w}'$, due to Lemma~\ref{le:C-density}.
In view of  the Mergelyan theorem (see, e.g., \cite[Ch.~III, Section~2]{Gaier87}),
every function in $C$ may uniformly be approximated by functions in
$\cR$. Thus, $\cR$ is dense in $X_w$ and in $X_{1/w}'$.
If $X$ is reflexive, we need first apply Corollary~\ref{co:weighted-reflexivity}
and then repeat the above arguments.
\end{proof}
%%%%%%%%%%%%%%%%%%%%%%%%%%%%%%%%%%%%%%%%%%%%%%%%%%%%%%%%%%%%%%%%%%%%%%%%%%%
\section{Analogs of the Muckenhoupt class}\label{sec:AX}
\subsection{Definitions}
Let $X$ be a Banach function space.
Fix $t\in\Gamma$. For a weight $w:\Gamma\to[0,\infty]$, put
\[
B_{t,R}(w):=\frac{1}{R}\|w\chi_{\Gamma(t,R)}\|_X\|\chi_{\Gamma(t,R)}/w\|_{X'}
\quad (R>0),
\]
where $\chi_{\Gamma(t,R)}$ is the characteristic function of the
portion $\Gamma(t,R)$. Consider the following classes of weights:
\[
A_X(\Gamma,t):=\Big\{w: \ \sup_{R>0}B_{t,R}(w)<\infty\Big\},
\quad
A_X(\Gamma):=\Big\{w: \ \sup_{t\in\Gamma}\sup_{R>0}B_{t,R}(w)<\infty\Big\}.
\]
Obviously, $A_X(\Gamma)\subset A_X(\Gamma,t)$ for $t\in\Gamma$. If
$X$ is a Lebesgue space $L^p,p\in(1,\infty)$, then $A_X(\Gamma)$ is the
Muckenhoupt class $A_p(\Gamma)$.
For a detailed discussion of Muckenhoupt weights on curves, see, e.g.,
\cite{BK97}. The classes $A_X(\Gamma,t)$ and $A_X(\Gamma)$ were defined
in \cite{K98} (see also \cite{K96,K00}) for rearrangement-invariant
spaces $X$. Here we assume only that $X$ is a Banach function space.
Ours definition is similar to a definition in \cite{Berezhnoi99}.
For others generalizations  (different from our) of the Muckenhoupt class
$A_p(\Gamma)$ in the setting of Orlicz and Lorentz spaces, see, e.g.,
\cite{GGKK98,KK91} and in the setting of Banach function spaces, see
\cite{KOPR93}.

With the help of H\"older's inequality (see Lemma~\ref{le:Hoelder}),
it is easy to show that $w\in A_X(\Gamma,t)$ implies
%%%
\begin{equation}\label{eq:loc-Carleson}
C_{\Gamma,t}:=\sup_{R>0}\frac{|\Gamma(t,R)|}{R}<\infty.
\end{equation}
%%%
We say that a rectifiable Jordan curve $\Gamma$ is \textit{locally a
Carleson curve at the point} $t\in\Gamma$ if (\ref{eq:loc-Carleson})
is satisfied. In that case the constant $C_{\Gamma,t}$ is referred to as the
\textit{local Carleson constant at the point} $t\in\Gamma$.
Analogously, if $w\in A_X(\Gamma)$, then
\[
C_\Gamma=\sup_{t\in\Gamma}C_{\Gamma,t}<\infty,
\]
that is, $\Gamma$ is a Carleson curve.
%%%%%%%%%%%%%%%%%%%%%%%%%%%%%%%%%%%%%%%%%%%%%%%%%%%%%%%%%%%%%%%%%%%%%%%%%%%
\subsection{Bounded and vanishing mean oscillation}
Let $\Gamma$ be a rectifiable Jordan curve. Let $f:\Gamma\to[-\infty,\infty]$
and $f\in L^1(\Gamma)$. Suppose $t\in\Gamma,\delta\in(0,\infty]$, and
$R\in(0,\infty)$. Put
%%%
\begin{eqnarray*}
\Omega_t(f,R) &:=& \frac{1}{|\Gamma(t,R)|}\int_{\Gamma(t,R)}f(\tau)|d\tau|,
\\
M_{\delta,t}(f) &:=& \sup_{0<R<\delta}
\frac{1}{|\Gamma(t,R)|}\int_{\Gamma(t,R)}|f(\tau)-\Omega_t(f,R)||d\tau|.
\end{eqnarray*}
%%%
A function $f$ is said to be of \textit{bounded mean oscillation
at the point} $t\in\Gamma$ if $\|f\|_{*,t}:=M_{\infty,t}(f)<\infty$.
In this case we will write $f\in BMO(\Gamma,t)$. A function $f\in BMO(\Gamma,t)$
has \textit{vanishing mean oscillation at the point} $t\in\Gamma$ if
\[
\lim_{\delta\to 0}M_{\delta,t}(f)=0.
\]
In that case we will write $f\in VMO(\Gamma,t)$.

One says that a function
$f:\Gamma\to[-\infty,\infty]$ is of \textit{bounded mean oscillation on}
$\Gamma$ if $f\in BMO(\Gamma,t)$ for all $t\in\Gamma$ and
\[
\|f\|_*:=\sup_{t\in\Gamma}\|f\|_{*,t}<\infty.
\]
The class of functions of bounded mean oscillation on $\Gamma$ is denoted by
$BMO(\Gamma)$. A function $f\in BMO(\Gamma)$ is said to be of
\textit{vanishing mean oscillation on} $\Gamma$ if
\[
\lim_{\delta\to 0}\sup_{t\in\Gamma}M_{\delta,t}(f)=0.
\]
The class of functions of vanishing mean oscillation on $\Gamma$ is denoted by
$VMO(\Gamma)$. Clearly, $BMO(\Gamma)\subset BMO(\Gamma,t)$ and
$VMO(\Gamma)\subset VMO(\Gamma,t)$ for every $t\in\Gamma$.
%%%%%%%%%%%%%%%%%%%%%%%%%%%%%%%%%%%%%%%%%%%%%%%%%%%%%%%%%%%%%%%%%%%%%%%%%%%%
\subsection{Bounded mean oscillation of logarithms of weights}
Let
\[
d_t:=\max_{\tau\in\Gamma}|\tau-t|.
\]
For a weight $w:\Gamma\to[0,\infty]$ such that $w\in X(\Gamma)$ and
$1/w\in X'(\Gamma)$, we have $w,1/w\in L^1(\Gamma)$. Then, taking into
account the obvious inequality $|\log x|\le x+1/x$ for $x\in (0,\infty)$,
we deduce that $\log w\in L^1$. For $t\in\Gamma$ and $R>0$, put
%%%
\begin{eqnarray*}
C(w,t,R) &:=& \exp(-\Omega_t(\log w,R))
\frac{\|w\chi_{\Gamma(t,R)}\|_{X}\|\chi_{\Gamma(t,R)}\|_{X'}}{|\Gamma(t,R)|},
\\
C'(w,t,R) &:=& \exp(\Omega_t(\log w,R))
\frac{\|\chi_{\Gamma(t,R)}\|_{X}\|\chi_{\Gamma(t,R)}/w\|_{X'}}{|\Gamma(t,R)|}.
\end{eqnarray*}
%%%
Clearly, these quantities are well defined.
%%%%%%%%%%%%%%%%%%%%%%%%%%%%%%%%%%%%%%%%%%%%%%%%%%%%%%%%%%%%%%%%%%%%%%%%%%%
\begin{lemma}\label{le:BMO-aux}
{\rm (a)} If $w\in A_X(\Gamma,t)$ and $1\in A_X(\Gamma,t)$, then
%%%
\begin{equation}\label{eq:BMO-aux-1}
1\le\sup_{R>0} C(w,t,R)<\infty,
\quad
1\le\sup_{R>0} C'(w,t,R)<\infty.
\end{equation}
%%%
{\rm (b)} If $w\in A_X(\Gamma)$ and $1\in A_X(\Gamma)$, then
%%%
\begin{equation}\label{eq:BMO-aux-2}
1\le\sup_{t\in\Gamma}\sup_{R>0} C(w,t,R)<\infty,
\quad
1\le\sup_{t\in\Gamma}\sup_{R>0} C'(w,t,R)<\infty.
\end{equation}
\end{lemma}
%%%%%%%%%%%%%%%%%%%%%%%%%%%%%%%%%%%%%%%%%%%%%%%%%%%%%%%%%%%%%%%%%%%%%%%%%%%
\begin{proof}
The proof is developed by similarity to \cite[Lemma~1.5]{K00}. Applying
Jensen's inequality (see, e.g., \cite[p.~78]{KR58}) and
H\"older's inequality (see Lemma~\ref{le:Hoelder}), we obtain
\[
\exp(\Omega_t(\log w,R))
\le
\frac{1}{|\Gamma(t,R)|}\int_{\Gamma(t,R)}
w(\tau)|d\tau|
\le
\frac{\|w\chi_{\Gamma(t,R)}\|_{X}\|\chi_{\Gamma(t,R)}\|_{X'}}{|\Gamma(t,R)|}.
\]
Hence,
%%%
\begin{equation}\label{eq:BMO-aux-3}
1\le C(w,t,R), \quad R>0.
\end{equation}
%%%
Analogously,
%%%
\begin{equation}\label{eq:BMO-aux-4}
1\le C'(w,t,R), \quad R>0.
\end{equation}
%%%
Inequalities (\ref{eq:BMO-aux-3}) and (\ref{eq:BMO-aux-4}) imply that
(\ref{eq:BMO-aux-1}) is equivalent to
%%%
\begin{equation}\label{eq:BMO-aux-5}
\sup_{R>0}\Big(C(w,t,R)C'(w,t,R)\Big)<\infty
\end{equation}
%%%
and (\ref{eq:BMO-aux-2}) is equivalent to
%%%
\begin{equation}\label{eq:BMO-aux-6}
\sup_{t\in\Gamma}\sup_{R>0}\Big(C(w,t,R)C'(w,t,R)\Big)<\infty.
\end{equation}
%%%
Since $\Gamma(t,R)=\Gamma$ for $R>d_t$, we have for every $t\in\Gamma$,
%%%
\begin{eqnarray}
\label{eq:BMO-aux-7}
\sup_{R>0} B_{t,R}(w)=\sup_{R\in(0,2d_t]} B_{t,R}(w),
\quad
\sup_{R>0} B_{t,R}(1)=\sup_{R\in(0,2d_t]} B_{t,R}(1),
\\
\label{eq:BMO-aux-8}
\sup_{R>0} \Big(C(w,t,R)C'(w,t,R)\Big)=\sup_{0<R\le 2d_t}\Big(C(w,t,R)C'(w,t,R)\Big).
\end{eqnarray}
%%%
Evidently, $R/2\le |\Gamma(t,R)|$ for $R\in(0,2d_t]$.
Taking into account the latter inequality and the definitions of $C(w,t,R), C'(w,t,R)$,
we get for $t\in\Gamma$ and $R\in(0,2d_t]$,
%%%
\begin{eqnarray*}
C(w,t,R)C'(w,t,R)
&\le&
\frac{\|w\chi_{\Gamma(t,R)}\|_X\|\chi_{\Gamma(t,R)}/w\|_{X'}}{|\Gamma(t,R)|}
\cdot
\frac{\|\chi_{\Gamma(t,R)}\|_X\|\chi_{\Gamma(t,R)}\|_{X'}}{|\Gamma(t,R)|}
\\
&\le&
4B_{t,R}(w)B_{t,R}(1).
\end{eqnarray*}
%%%
Therefore,
%%%
\begin{eqnarray}
&&
\sup_{R\in(0,2d_t]}\Big(C(w,t,R)C'(w,t,R)\Big)
\le
4\Big(\sup_{R\in(0,2d_t]} B_{t,R}(w)\Big)
 \Big(\sup_{R\in(0,2d_t]} B_{t,R}(1)\Big).
\label{eq:BMO-aux-10}
\end{eqnarray}
%%%
From (\ref{eq:BMO-aux-7})--(\ref{eq:BMO-aux-10})
it follows that
%%%
\begin{eqnarray}
&&
\sup_{R>0}\Big(C(w,t,R)C'(w,t,R)\Big)
\le
4\Big(\sup_{R>0} B_{t,R}(w)\Big)
 \Big(\sup_{R>0} B_{t,R}(1)\Big),
\label{eq:BMO-aux-11}
\\
&&
\sup_{t\in\Gamma}\sup_{R>0}\Big(C(w,t,R)C'(w,t,R)\Big)
\le
4\Big(\sup_{t\in\Gamma}\sup_{R>0} B_{t,R}(w)\Big)
 \Big(\sup_{t\in\Gamma}\sup_{R>0} B_{t,R}(1)\Big).
\label{eq:BMO-aux-12}
\end{eqnarray}
%%%

(a) If $w\in A_X(\Gamma,t)$ and $1\in A_X(\Gamma,t)$, then (\ref{eq:BMO-aux-11})
implies (\ref{eq:BMO-aux-5}), but we have shown that (\ref{eq:BMO-aux-5})
is equivalent to (\ref{eq:BMO-aux-1}). Part (a) is proved.
Part (b) is proved similarly by using (\ref{eq:BMO-aux-12})
and the equivalence of (\ref{eq:BMO-aux-6}) and (\ref{eq:BMO-aux-2}).
\end{proof}
%%%%%%%%%%%%%%%%%%%%%%%%%%%%%%%%%%%%%%%%%%%%%%%%%%%%%%%%%%%%%%%%%%%%%%%%%%%
\begin{lemma}\label{le:BMO}
{\rm (a)} If $w\in A_X(\Gamma,t)$ and $1\in A_X(\Gamma,t)$, then
$\log w\in BMO(\Gamma,t)$.

\noindent
{\rm (b)} If $w\in A_X(\Gamma)$ and $1\in A_X(\Gamma)$, then
$\log w\in BMO(\Gamma)$.
\end{lemma}
%%%%%%%%%%%%%%%%%%%%%%%%%%%%%%%%%%%%%%%%%%%%%%%%%%%%%%%%%%%%%%%%%%%%%%%%%%%
\begin{proof}
This statement is proved by analogy with \cite[Lemma~1.6]{K00} (see also
\cite[Proposition~2.4]{BK97}). Put $\Omega_t(R):=\Omega_t(\log w,R)$,
%%%
\begin{eqnarray*}
\Gamma_+(t,R)&:=&
\Big\{\tau\in\Gamma(t,R):\log w(\tau)\ge\Omega_t(R)\Big\},
\\
\Gamma_-(t,R)&:=&
\Big\{\tau\in\Gamma(t,R):\log w(\tau)<\Omega_t(R)\Big\}.
\end{eqnarray*}
%%%
Due to Jensen's inequality \cite[p.~78]{KR58},
%%%
\begin{eqnarray}
\label{eq:BMO-1}
&&
\exp\left(
\frac{1}{|\Gamma(t,R)|}\int_{\Gamma(t,R)}|\log w(\tau)-\Omega_t(R)||d\tau|
\right)
\\
&&
\le
\frac{1}{|\Gamma(t,R)|}\int_{\Gamma_+(t,R)}
\exp\Big(\log w(\tau)-\Omega_t(R)\Big)|d\tau|
\nonumber\\
&&
+
\frac{1}{|\Gamma(t,R)|}\int_{\Gamma_-(t,R)}
\exp\Big(-(\log w(\tau)-\Omega_t(R))\Big)|d\tau|
\nonumber\\
&&
\le
\frac{1}{|\Gamma(t,R)|}\int_\Gamma
\exp\Big(\log w(\tau)-\Omega_t(R)\Big)\chi_{\Gamma(t,R)}(\tau)|d\tau|
\nonumber\\
&&
+
\frac{1}{|\Gamma(t,R)|}\int_\Gamma
\exp\Big(-(\log w(\tau)-\Omega_t(R))\Big)\chi_{\Gamma(t,R)}(\tau)|d\tau|.
\nonumber
\end{eqnarray}
%%%
Applying H\"older's inequality (see Lemma~\ref{le:Hoelder}) to the first
term on the right of (\ref{eq:BMO-1}), we get
%%%
\begin{eqnarray}
\label{eq:BMO-2}
&&
\frac{1}{|\Gamma(t,R)|}\int_\Gamma
\exp\Big(\log w(\tau)-\Omega_t(R)\Big)\chi_{\Gamma(t,R)}(\tau)|d\tau|
\\
&&
\le
\Big\|\exp\Big(\log w(\cdot)-\Omega_t(R)\Big)\chi_{\Gamma(t,R)}(\cdot)\Big\|_X
\frac{\|\chi_{\Gamma(t,R)}\|_{X'}}{|\Gamma(t,R)|}
\nonumber\\
&&
=e^{-\Omega_t(R)}
\frac{\|\chi_{\Gamma(t,R)}\|_X\|\chi_{\Gamma(t,R)}\|_{X'}}{|\Gamma(t,R)|}
=C(w,t,R).
\nonumber
\end{eqnarray}
%%%
Analogously,
%%%
\begin{equation}\label{eq:BMO-3}
\frac{1}{|\Gamma(t,R)|}\int_\Gamma
\exp\Big(-(\log w(\tau)-\Omega_t(R))\Big)\chi_{\Gamma(t,R)}(\tau)|d\tau|
\le C'(w,t,R).
\end{equation}
%%%
Combining (\ref{eq:BMO-1})--(\ref{eq:BMO-3}), we see that for every
$t\in\Gamma$ and $R>0$,
\[
\exp\left(
\frac{1}{|\Gamma(t,R)|}\int_{\Gamma(t,R)}|\log w(\tau)-\Omega_t(R)||d\tau|
\right)
\le
C(w,t,R)+C'(w,t,R).
\]
Consequently,
%%%
\begin{eqnarray}
&&
\|\log w\|_{*,t}\le\log\left(\sup_{R>0}C(w,t,R)+\sup_{R>0}C'(w,t,R)\right),
\quad t\in\Gamma,
\label{eq:BMO-4}
\\
&&
\|\log w\|_*\le\log\left(\sup_{t\in\Gamma}\sup_{R>0}C(w,t,R)+
\sup_{t\in\Gamma}\sup_{R>0}C'(w,t,R)\right).
\label{eq:BMO-5}
\end{eqnarray}
%%%
Statement (a) follows from Lemma~\ref{le:BMO-aux}(a) and (\ref{eq:BMO-4}).
Statement (b) follows from Lemma~\ref{le:BMO-aux}(b) and (\ref{eq:BMO-5}).
\end{proof}
%%%%%%%%%%%%%%%%%%%%%%%%%%%%%%%%%%%%%%%%%%%%%%%%%%%%%%%%%%%%%%%%%%%%%%%%%%%
For rearrangement-invariant Banach function spaces $X$, by using
(\ref{eq:RI-fundamental}), we infer that $w\in A_X(\Gamma)$
implies $1\in A_X(\Gamma)$. In that case, by Lemma~\ref{le:BMO}(b),
if $w\in A_X(\Gamma)$, then $\log w\in BMO(\Gamma)$.
This result was obtained in \cite[Lemma~1.6]{K00}. Note that
for Lebesgue spaces $L^p, 1<p<\infty$, and Muckenhoupt classes $A_p(\Gamma)$
this fact is well known (see, e.g., \cite[Proposition~2.4]{BK97}).
%%%%%%%%%%%%%%%%%%%%%%%%%%%%%%%%%%%%%%%%%%%%%%%%%%%%%%%%%%%%%%%%%%%%%%%%%%%
\section{Indices of submultiplicative functions associated \\
with weighted Banach function spaces}\label{sec:indices}
\subsection{Submultiplicative functions and their indices}
Following \cite[Section~1.4]{BK97}, we say a function
$\Phi:(0,\infty)\to(0,\infty]$ is \textit{regular} if it is bounded
in an open neighborhood of $1$. A function
$\Phi:(0,\infty)\to(0,\infty]$ is said to be \textit{submultiplicative} if
\[
\Phi(xy)\le\Phi(x)\Phi(y)
\quad\mbox{for all}\quad x,y\in(0,\infty).
\]
It is easy to show that if $\Phi$ is regular and submultiplicative, then
$\Phi$ is bounded away from zero in some open neighborhood of $1$.
Moreover, in this case $\Phi(x)$ is finite for all $x\in(0,\infty)$.
Given a regular and submultiplicative function $\Phi:(0,\infty)\to(0,\infty)$,
one defines
\[
\alpha(\Phi):=\sup_{x\in(0,1)}\frac{\log\Phi(x)}{\log x},
\quad
\beta(\Phi):=\inf_{x\in(1,\infty)}\frac{\log\Phi(x)}{\log x}.
\]
Clearly, $-\infty<\alpha(\Phi)$ and $\beta(\Phi)<\infty$.
%%%%%%%%%%%%%%%%%%%%%%%%%%%%%%%%%%%%%%%%%%%%%%%%%%%%%%%%%%%%%%%%%%%%%%%%%%%
\begin{theorem}\label{th:submult}
{\rm (see \cite[Theorem~1.13]{BK97})}. % or \cite[Ch.~2, Theorem~1.3]{KPS78}).}
If $\Phi:(0,\infty)\to(0,\infty)$ is regular and submultiplicative, then
\[
\alpha(\Phi)=\lim_{x\to 0}\frac{\log\Phi(x)}{\log x},
\quad
\beta(\Phi)=\lim_{x\to\infty}\frac{\log\Phi(x)}{\log x}
\]
and $-\infty<\alpha(\Phi)\le\beta(\Phi)<+\infty$.
\end{theorem}
%%%%%%%%%%%%%%%%%%%%%%%%%%%%%%%%%%%%%%%%%%%%%%%%%%%%%%%%%%%%%%%%%%%%%%%%%%
The quantities $\alpha(\Phi)$ and $\beta(\Phi)$ are called the \textit{lower}
and  \textit{upper indices of the regular and submultiplicative function}
$\Phi$, respectively.
%%%%%%%%%%%%%%%%%%%%%%%%%%%%%%%%%%%%%%%%%%%%%%%%%%%%%%%%%%%%%%%%%%%%%%%%%%%
\subsection{Spirality indices}
In this subsection we mainly follow \cite[Ch.~1]{BK97}. Fix $t\in\Gamma$.
Suppose $\psi:\Gamma\setminus\{t\}\to(0,\infty)$ is a continuous function. Put
\[
F_\psi(R_1,R_2):=\max_{\tau\in\Gamma,|\tau-t|=R_1}\psi(\tau)\Big/
\min_{\tau\in\Gamma,|\tau-t|=R_2}\psi(\tau),
\quad
R_1,R_2\in(0,d_t].
\]
By \cite[Lemma~1.15]{BK97}, the function
\[
(W_t\psi)(x):=\left\{\begin{array}{ll}
\displaystyle
\sup_{0<R\le d_t}F_\psi(xR,R), & x\in(0,1],\\
\displaystyle
\sup_{0<R\le d_t}F_\psi(R,x^{-1}R), & x\in(1,\infty).
\end{array}
\right.
\]
is submultiplicative. For $t\in\Gamma$, we have,
\[
\tau- t=|\tau-t|e^{i\arg(\tau-t)},
\quad
\tau\in\Gamma\setminus\{t\},
\]
and the argument $\arg(\tau-t)$ may be chosen to be a continuous function
of $\tau\in\Gamma\setminus\{t\}$. Consider
\[
\eta_t(\tau):=e^{-\arg(\tau-t)}.
\]
Using the local Carleson constant $C_{\Gamma,t}$ instead of the global
Carleson constant $C_\Gamma$, we can obtain the following local
versions of \cite[Theorem~1.10 and Lemma~1.17]{BK97}.
%%%%%%%%%%%%%%%%%%%%%%%%%%%%%%%%%%%%%%%%%%%%%%%%%%%%%%%%%%%%%%%%%%%%%%%%%%%
\begin{lemma}\label{le:Seifullaev}
If $\Gamma$ is locally a Carleson curve at $t\in\Gamma$, then
%%%
\[
\arg(\tau-t)=O(-\log|\tau-t|)
\quad\mbox{as}\quad\tau\to t.
\]
\end{lemma}
%%%%%%%%%%%%%%%%%%%%%%%%%%%%%%%%%%%%%%%%%%%%%%%%%%%%%%%%%%%%%%%%%%%%%%%%%%%
\begin{lemma}\label{le:W-reg}
If $\Gamma$ is locally a Carleson curve at $t\in\Gamma$, then the
submultiplicative function $W_t\eta_t$ is regular.
\end{lemma}
%%%%%%%%%%%%%%%%%%%%%%%%%%%%%%%%%%%%%%%%%%%%%%%%%%%%%%%%%%%%%%%%%%%%%%%%%%%
Under the assumptions of Lemma~\ref{le:W-reg}, by Theorem~\ref{th:submult},
there exist the \textit{spirality indices}
\[
\delta_t^-:=\alpha(W_t\eta_t),
\quad
\delta_t^+:=\beta(W_t\eta_t)
\]
\textit{of the curve} $\Gamma$ \textit{at the point} $t$ (see \cite[Ch.~1]{BK97}).
If, in addition,
%%%
\[
\arg(\tau-t)=-\delta_t\log|\tau-t|+O(1)
\quad\mbox{as}\quad\tau\to t,
\]
where $\delta_t\in\R$, then $\delta_t^-=\delta_t^+=\delta_t$
(see \cite[Section~1.6]{BK97}). Examples of Carleson curves with distinct
spirality indices are also given there.

On a rectifiable Jordan curve we have $d\tau=e^{i\theta_\Gamma(\tau)}|d\tau|$
where $\theta_\Gamma(\tau)$ is the angle between the positively oriented
real axis and the naturally oriented tangent of $\Gamma$ at $\tau$
(which exists almost everywhere). A rectifiable Jordan curve $\Gamma$
is said to be a \textit{Lyapunov curve} if
\[
|\theta_\Gamma(\tau)-\theta_\Gamma(t)|\le c|\tau-t|^\mu
\]
for some constants $c>0, \mu\in(0,1)$  and all $\tau,t\in\Gamma$.
If $\theta_\Gamma$ is a function of bounded variation on
$\Gamma$, then the curve $\Gamma$ is called a \textit{Radon curve}
(or a \textit{curve of bounded rotation}). It is very well known
that Lyapunov curves are smooth, but Radon curves may have
at most countable set of corner points (or even cusps).
All Lyapunov curves and Radon curves without cusps are
Carleson curves (see, e.g., \cite[Section~2.3]{Khvedelidze75}).
The next statement is well known.
%%%%%%%%%%%%%%%%%%%%%%%%%%%%%%%%%%%%%%%%%%%%%%%%%%%%%%%%%%%%%%%%%%%%%%%%%%%
\begin{proposition}\label{pr:nice-curves}
If $\Gamma$ is either a Lyapunov Jordan curve or a Radon Jordan curve,
then for every $t\in\Gamma$,
\[
\arg(\tau-t)=O(1)
\quad\mbox{as}\quad\tau\to t,
\]
and, therefore, $\delta_t^-=\delta_t^+=0$.
\end{proposition}
%%%%%%%%%%%%%%%%%%%%%%%%%%%%%%%%%%%%%%%%%%%%%%%%%%%%%%%%%%%%%%%%%%%%%%%%%%%
\subsection{Indices of powerlikeness}
To investigate whether the weight $|(\tau-t)^\gamma|w(\tau)$ with arbitrary
$\gamma\in\C$ belongs to the Muckenhoupt class $A_p(\Gamma)$,
A.~B\"ottcher and Yu.~I.~Karlovich introduced submultiplicative functions
$V_tw$ and $V_t^0w$  associated with local properties of the
weight $w$ at the point $t\in\Gamma$ (see \cite[Ch.~3]{BK97}).

Let $w$ be a weight on $\Gamma$ such that $\log w\in L^1(\Gamma(t,R))$ for
every $R\in(0,d_t]$. Put
\[
H_w(R_1,R_2):=\exp(\Omega_t(\log w,R_1))/\exp(\Omega_t(\log w,R_2)),
\quad R_1,R_2\in(0,d_t].
\]
Consider the functions
%%%
\begin{eqnarray*}
(V_tw)(x) &:=&
\left\{\begin{array}{ll}
\displaystyle
\sup_{0<R\le d_t} H_w(xR,R), & x\in(0,1],\\
\displaystyle
\sup_{0<R\le d_t} H_w(R,x^{-1}R), & x\in(1,\infty),
\end{array}\right.
\\
(V_t^0w)(x)&:=&\limsup_{R\to 0} H_w(xR,R),
\quad x\in(0,\infty).
\end{eqnarray*}
%%%%%%%%%%%%%%%%%%%%%%%%%%%%%%%%%%%%%%%%%%%%%%%%%%%%%%%%%%%%%%%%%%%%%%%%%%%
\begin{lemma}\label{le:V-sub}
The function $V_tw$ is submultiplicative. If $V_tw$ is regular, then
$V_t^0w$ is regular and submultiplicative. Moreover,
$\alpha(V_t^0w)=\alpha(V_tw)$ and $\beta(V_t^0w)=\beta(V_tw)$.
\end{lemma}
%%%%%%%%%%%%%%%%%%%%%%%%%%%%%%%%%%%%%%%%%%%%%%%%%%%%%%%%%%%%%%%%%%%%%%%%%%%
\begin{lemma}\label{le:V-reg}
If $\Gamma$ is locally a Carleson curve at $t\in\Gamma$ and
$\log w\in BMO(\Gamma,t)$, then $V_tw$ and $V_t^0w$ are regular.
\end{lemma}
%%%%%%%%%%%%%%%%%%%%%%%%%%%%%%%%%%%%%%%%%%%%%%%%%%%%%%%%%%%%%%%%%%%%%%%%%%%
Lemmas~\ref{le:V-sub} and~\ref{le:V-reg} are proved by analogy
with \cite[Lemma~3.5(a)]{BK97} and \cite[Lemma~3.2(a)]{BK97}.
These statements are stated in \cite{BK97} under the assumption that
$\Gamma$ is a Carleson curve. But Lemma~\ref{le:V-sub} is valid
for arbitrary rectifiable curves $\Gamma$. Since Lemma~\ref{le:V-reg}
has a ``local nature'', we may use the ``local'' Carleson constant
$C_{\Gamma,t}$ instead of the ``global'' Carleson constant $C_\Gamma$
in its proof. Under the assumptions of Lemma~\ref{le:V-reg},
in view of Theorem~\ref{th:submult}, for the weight $w$, there exist
the \textit{indices of powerlikeness}
%%%
\begin{equation}\label{eq:ind_powerlikeness}
\mu_t:=\alpha(V_t^0w)=\alpha(V_tw),
\quad
\nu_t:=\beta(V_t^0w)=\beta(V_tw)
\end{equation}
%%%
\textit{at the point} $t\in\Gamma$.

Obviously, for a power weight $w(\tau)=|\tau-t|^{\lambda_t}$, the indices
of powerlikeness equal $\mu_t=\nu_t=\lambda_t$. Nontrivial examples
of weights with distinct indices of powerlikeness are given in
\cite[Examples~3.24--3.28]{BK97}.
%%%%%%%%%%%%%%%%%%%%%%%%%%%%%%%%%%%%%%%%%%%%%%%%%%%%%%%%%%%%%%%%%%%%%%%%%%%
\begin{lemma}\label{le:VMO-indices}
{\rm (see \cite[Lemma~2.4]{K00}).}
If $\Gamma$ is locally a Carleson curve at $t\in\Gamma$ and
$\log w\in VMO(\Gamma,t)$, then $\mu_t=\nu_t=0$.
\end{lemma}
%%%%%%%%%%%%%%%%%%%%%%%%%%%%%%%%%%%%%%%%%%%%%%%%%%%%%%%%%%%%%%%%%%%%%%%%%%%
\subsection{Submultiplicative functions associated with weighted Banach
function spaces}
Let $\Gamma$ be a  rectifiable Jordan curve and let $X$ be a Banach function
space. Fix $t\in\Gamma$ and consider the portion of the curve $\Gamma$
in the annulus
\[
\Delta(t,R):=\Gamma(t,R)\setminus\Gamma(t,R/2),\quad R>0.
\]
Clearly,
%%%
\begin{equation}\label{eq:Delta-lower}
R/2\le |\Delta(t,R)|,\quad R\in(0,d_t].
\end{equation}
%%%
On the other hand, if $\Gamma$ is locally a Carleson curve at $t\in\Gamma$, then
%%%
\begin{equation}\label{eq:Delta-upper}
|\Delta(t,R)|\le|\Gamma(t,R)|\le C_{\Gamma,t}R,
\quad R>0.
\end{equation}
%%%
Suppose $w:\Gamma\to[0,\infty]$ is a weight such that $w\chi_{\Delta(t,R)}\in X$
and $\chi_{\Delta(t,R)}/w\in X'$ for all $R\in(0,d_t]$. We denote
\[
G_w(R_1,R_2):=
\frac{\|w\chi_{\Delta(t,R_1)}\|_X\|\chi_{\Delta(t,R_2)}/w\|_{X'}}{|\Delta(t,R_2)|},
\quad R_1,R_2\in(0,d_t].
\]
Define the following functions (see \cite[Section~5]{K98}):
%%%
\begin{eqnarray*}
(Q_tw)(x) &:=&
\left\{\begin{array}{ll}
\displaystyle
\sup_{0<R\le d_t} G_w(xR,R), & x\in(0,1],\\
\displaystyle
\sup_{0<R\le d_t} G_w(R,x^{-1}R), & x\in(1,\infty),
\end{array}\right.
\\
(Q_t^0w)(x)&:=&\limsup_{R\to 0} G_w(xR,R),
\quad x\in(0,\infty).
\end{eqnarray*}
%%%%%%%%%%%%%%%%%%%%%%%%%%%%%%%%%%%%%%%%%%%%%%%%%%%%%%%%%%%%%%%%%%%%%%%%%%%
\begin{lemma}\label{le:Q-sub}
The function $Q_tw$ is submultiplicative. If $Q_tw$ is regular, then
$Q_t^0w$ is regular and submultiplicative. Moreover,
$\alpha(Q_t^0w)=\alpha(Q_tw),\beta(Q_t^0w)=\beta(Q_tw)$.
\end{lemma}
%%%%%%%%%%%%%%%%%%%%%%%%%%%%%%%%%%%%%%%%%%%%%%%%%%%%%%%%%%%%%%%%%%%%%%%%%%%
\begin{lemma}\label{le:Q-reg}
If $w\in A_X(\Gamma,t)$, then $Q_tw$ and $Q_t^0w$ are regular. Moreover,
\[
0\le\alpha(Q_tw)=\alpha(Q_t^0w)\le\beta(Q_t^0w)=\beta(Q_tw)\le 1.
\]
\end{lemma}
%%%%%%%%%%%%%%%%%%%%%%%%%%%%%%%%%%%%%%%%%%%%%%%%%%%%%%%%%%%%%%%%%%%%%%%%%%%
These statements are proved in \cite[Lemmas~5.1--5.2]{K98}
and \cite[Theorem~5.3]{K98}, respectively, under the assumption that
$X$ is rearrangement-invariant. But, actually, we did not use this
assumption in those proofs. So we can literally repeat the proofs for
arbitrary Banach function spaces.
%%%%%%%%%%%%%%%%%%%%%%%%%%%%%%%%%%%%%%%%%%%%%%%%%%%%%%%%%%%%%%%%%%%%%%%%%%%
\section{Relations between indices}\label{sec:relations}
\subsection{Case of general Banach function spaces}
Let $\Gamma$ be a rectifiable Jordan curve, let $X$ be a Banach function
space, and let $t\in\Gamma$.
%%%%%%%%%%%%%%%%%%%%%%%%%%%%%%%%%%%%%%%%%%%%%%%%%%%%%%%%%%%%%%%%%%%%%%%%%%%
\begin{theorem}\label{th:VW}
Suppose $w:\Gamma\to[0,\infty]$ is a weight such that $\log w\in L^1(\Gamma(t,R))$
for every $R\in(0,d_t]$ and
$\psi:\Gamma\setminus\{t\}\to(0,\infty)$ is a continuous function.
If the functions $V_tw$ and $W_t\psi$ are regular, then the function
$V_t(\psi w)$ is regular too. Moreover,
\[
\begin{split}
\alpha(V_tw)+\alpha(W_t\psi)
\le
\alpha(V_t(\psi w))
\le
\min\Big\{\alpha(V_tw)+\beta(W_t\psi),\beta(V_tw)+\alpha(W_t\psi)\Big\},
\\
\beta(V_tw)+\beta(W_t\psi)
\ge
\beta(V_t(\psi w))
\ge
\max\Big\{\alpha(V_tw)+\beta(W_t\psi),\beta(V_tw)+\alpha(W_t\psi)\Big\}.
\end{split}
\]
\end{theorem}

This statement is proved similarly to \cite[Lemma~3.17]{BK97}.
%%%%%%%%%%%%%%%%%%%%%%%%%%%%%%%%%%%%%%%%%%%%%%%%%%%%%%%%%%%%%%%%%%%%%%%%%%%
\begin{theorem}\label{th:QW}
Suppose $w:\Gamma\to[0,\infty]$ is a weight such that
$w\chi_{\Delta(t,R)}\in X$ and $\chi_{\Delta(t,R)}/w\in X'$ for every
$R\in(0,d_t]$ and $\psi:\Gamma\setminus\{t\}\to(0,\infty)$ is a continuous
function. If the functions $Q_tw$ and $W_t\psi$ are regular, then the function
$Q_t(\psi w)$ is regular too. Moreover,
\[
\begin{split}
\alpha(Q_tw)+\alpha(W_t\psi)
\le
\alpha(Q_t(\psi w))
\le
\min\Big\{\alpha(Q_tw)+\beta(W_t\psi),\beta(Q_tw)+\alpha(W_t\psi)\Big\},
\\
\beta(Q_tw)+\beta(W_t\psi)
\ge
\beta(Q_t(\psi w))
\ge
\max\Big\{\alpha(Q_tw)+\beta(W_t\psi),\beta(Q_tw)+\alpha(W_t\psi)\Big\}.
\end{split}
\]
\end{theorem}
%%%%%%%%%%%%%%%%%%%%%%%%%%%%%%%%%%%%%%%%%%%%%%%%%%%%%%%%%%%%%%%%%%%%%%%%%%%

This theorem is proved in \cite[Theorem~5.8]{K98} for rearrangement-invariant
Banach function spaces. The proof given there does not use the
rearrangement-invariant property of the space, so it works for an arbitrary
Banach function space.
%%%%%%%%%%%%%%%%%%%%%%%%%%%%%%%%%%%%%%%%%%%%%%%%%%%%%%%%%%%%%%%%%%%%%%%%%%%
\begin{lemma}\label{le:BoKa}
If $\Gamma$ is locally a Carleson curve at $t\in\Gamma$ and
$\log w\in BMO(\Gamma,t)$, then for every $R\in(0,d_t]$,
\[
\exp(\Omega_t(\log w,R))
\le
\frac{C_t}{|\Delta(t,R)|}\int_{\Delta(t,R)} w(\tau)|d\tau|
\]
where $C_t:=\exp(2C_{\Gamma,t}\|\log w\|_{*,t})<\infty$.
\end{lemma}
%%%%%%%%%%%%%%%%%%%%%%%%%%%%%%%%%%%%%%%%%%%%%%%%%%%%%%%%%%%%%%%%%%%%%%%%%%%
The proof is actually given in \cite[Lemma~3.2(b)]{BK97}.
%%%%%%%%%%%%%%%%%%%%%%%%%%%%%%%%%%%%%%%%%%%%%%%%%%%%%%%%%%%%%%%%%%%%%%%%%%%
\begin{theorem}\label{th:QV-lower}
Let $\Gamma$ be locally a Carleson curve at $t\in\Gamma$ and let
$w:\Gamma\to[0,\infty]$ be a weight such that
$w\chi_{\Delta(t,R)}\in X, \chi_{\Delta(t,R)}/w\in X'$ for every
$R\in(0,d_t]$ and $\log w\in BMO(\Gamma,t)$. If $Q_tw$ and
$Q_t1$ are regular, then
%%%
\begin{equation}\label{eq:QV-lower-1}
\alpha(Q_tw)\le\mu_t+\beta(Q_t1),
\quad
\nu_t+\alpha(Q_t1)\le\beta(Q_tw).
\end{equation}
\end{theorem}
%%%%%%%%%%%%%%%%%%%%%%%%%%%%%%%%%%%%%%%%%%%%%%%%%%%%%%%%%%%%%%%%%%%%%%%%%%%
\begin{proof}
The proof is developed by analogy with \cite[Theorem~2.6]{K00}.
From Lemma~\ref{le:BoKa} and H\"older's inequality (see Lemma~\ref{le:Hoelder})
we see that for every $R\in(0,d_t]$,
%%%
\begin{eqnarray}
\exp(\Omega_t(\log w,R))
&\le& C_t
\frac{\|w\chi_{\Delta(t,R)}\|_X\|\chi_{\Delta(t,R)}\|_{X'}}{|\Delta(t,R)|},
\label{eq:QV-lower-2}
\\
\exp(-\Omega_t(\log w,R))
&\le& C_t
\frac{\|\chi_{\Delta(t,R)}\|_X\|\chi_{\Delta(t,R)}/w\|_{X'}}{|\Delta(t,R)|}.
\label{eq:QV-lower-3}
\end{eqnarray}
%%%
From (\ref{eq:QV-lower-2}) and (\ref{eq:QV-lower-3}) it follows that
for $x\in(0,1]$ and $R\in(0,d_t]$,
%%%
\begin{eqnarray}
\label{eq:QV-lower-4}
&&
H_w(xR,R)
=
\exp(\Omega_t(\log w,xR))\exp(-\Omega_t(\log w,R))
\\
&&\le
C_t^2
\frac{\|w\chi_{\Delta(t,xR)}\|_X\|\chi_{\Delta(t,R)}/w\|_{X'}}{|\Delta(t,R)|}
\cdot
\frac{\|\chi_{\Delta(t,R)}\|_X\|\chi_{\Delta(t,xR)}\|_{X'}}{|\Delta(t,xR)|}
\nonumber\\
&&=
C_t^2G_w(xR,R)G_1(R,xR).
\nonumber
\end{eqnarray}
%%%
Then, taking the supremum over all $R\in(0,d_t]$, we obtain for $x\in(0,1]$,
%%%
\begin{equation}\label{eq:QV-lower-5}
(V_tw)(x)\le C_t(Q_tw)(x)(Q_t1)(x^{-1}).
\end{equation}
%%%
Analogously, for $x\in(1,\infty)$ and $R\in(0,d_t]$,
%%%
\begin{equation}\label{eq:QV-lower-6}
H_w(R,x^{-1}R)\le C_t^2G_w(R,x^{-1}R)G_1(x^{-1}R,R).
\end{equation}
%%%
Taking the supremum over all $R\in(0,d_t]$, we arrive at (\ref{eq:QV-lower-5})
for $x\in(1,\infty)$.
By Lemmas~\ref{le:V-sub}--\ref{le:V-reg}, the function $V_tw$ is regular
and submultiplicative. By Lemma~\ref{le:Q-sub}, the functions $Q_tw$
and $Q_t1$ are submultiplicative, they are regular, due to the assumption
of the theorem. Therefore, in view of Theorem~\ref{th:submult}, the indices
$\alpha(Q_tw), \beta(Q_tw)$; $\alpha(Q_t1), \beta(Q_t1)$; and
$\alpha(V_tw),\beta(V_tw)$ exist and are well defined.

From (\ref{eq:QV-lower-5}) it follows that
%%%
\begin{eqnarray*}
&&
\frac{\log (V_tw)(x)}{\log x}
\ge
\frac{\log C_t^2}{\log x}+\frac{\log (Q_tw)(x)}{\log x}
-
\frac{\log(Q_t1)(x^{-1})}{\log x^{-1}},
\quad x\in(0,1],
\\
&&
\frac{\log (V_tw)(x)}{\log x}
\le
\frac{\log C_t^2}{\log x}+\frac{\log (Q_tw)(x)}{\log x}
-
\frac{\log(Q_t1)(x^{-1})}{\log x^{-1}},
\quad x\in(1,\infty).
\end{eqnarray*}
%%%
Passing to the limit in the latter inequalities as $x\to 0$ (respectively,
as $x\to\infty$), we obtain, respectively,
\[
\mu_t=\alpha(V_tw)\ge\alpha(Q_tw)-\beta(Q_t1),
\quad
\nu_t=\beta(V_tw)\le\beta(Q_tw)-\alpha(Q_t1).
\]
So, we arrive at (\ref{eq:QV-lower-1}).
\end{proof}
%%%%%%%%%%%%%%%%%%%%%%%%%%%%%%%%%%%%%%%%%%%%%%%%%%%%%%%%%%%%%%%%%%%%%%%%%%%
\begin{theorem}\label{th:QV}
If $w\in A_X(\Gamma,t)$ and $1\in A_X(\Gamma,t)$, then
\begin{eqnarray}
&&
\alpha(Q_t1)+\mu_t\le\alpha(Q_tw)\le\min\Big\{\alpha(Q_t1)+\nu_t,\beta(Q_t1)+\mu_t\Big\},
\label{eq:QV-1}
\\
&&
\beta(Q_t1)+\nu_t\ge\beta(Q_tw)\ge\max\Big\{\alpha(Q_t1)+\nu_t,\beta(Q_t1)+\mu_t\Big\}.
\label{eq:QV-2}
\end{eqnarray}
\end{theorem}
%%%%%%%%%%%%%%%%%%%%%%%%%%%%%%%%%%%%%%%%%%%%%%%%%%%%%%%%%%%%%%%%%%%%%%%%%%%
\begin{proof}
The idea of the proof is borrowed from \cite[Theorems~2.6 and~2.7]{K00}.
From Lemmas~\ref{le:Q-sub}--\ref{le:Q-reg} it follows that the functions
$Q_tw$ and $Q_t1$ are regular and submultiplicative. On the other hand,
by Lemma~\ref{le:BMO}(a), $\log w\in BMO(\Gamma,t)$. Therefore, by
Lemmas~\ref{le:V-sub}--\ref{le:V-reg}, the function $V_tw$ is regular and
submultiplicative. Thus, all the indices
\[
\alpha(Q_t1),\quad
\beta(Q_t1),\quad
\alpha(Q_tw),\quad
\beta(Q_tw),\quad
\mu_t=\alpha(V_tw),\quad
\nu_t=\beta(V_tw)
\]
are well defined. By Theorem~\ref{th:QV-lower},
%%%
\begin{equation}\label{eq:QV-3}
\alpha(Q_tw)\le\mu_t+\beta(Q_t1),
\quad
\nu_t+\alpha(Q_t1)\le\beta(Q_tw).
\end{equation}
%%%

If $1\in A_X(\Gamma,t)$, then from the lattice property it follows that for every $R>0$,
%%%
\begin{eqnarray}
\label{eq:QV-4}
\frac{1}{R}\|\chi_{\Delta(t,R)}\|_X\|\chi_{\Delta(t,R)}\|_{X'}
&\le&
\frac{1}{R}\|\chi_{\Gamma(t,R)}\|_X\|\chi_{\Gamma(t,R)}\|_{X'}
\\
&\le&
\sup_{R>0}B_{t,R}(1)=:B_t(1).
\nonumber
\end{eqnarray}
%%%
Combining (\ref{eq:QV-4}) and (\ref{eq:Delta-lower}), we arrive at
\[
\|\chi_{\Delta(t,R)}\|_X\|\chi_{\Delta(t,R)}\|_{X'}\le 2B_t(1)|\Delta(t,R)|,
\quad R\in(0,d_t].
\]
Then we have for $x\in(0,1]$,
%%%
\begin{eqnarray}
\label{eq:QV-6}
\frac{1}{G_1(R,xR)}
&=&
\frac{|\Delta(t,xR)|}{\|\chi_{\Delta(t,R)}\|_X\|\chi_{\Delta(t,xR)}\|_{X'}}
\ge
\frac{\|\chi_{\Delta(t,xR)}\|_X\|\chi_{\Delta(t,R)}\|_{X'}}{(2B_t(1))^2|\Delta(t,R)|}
\\
&=&
(2B_t(1))^{-2}G_1(xR,R).
\nonumber
\end{eqnarray}
%%%
Analogously, we deduce that for $x\in(1,\infty)$,
%%%
\begin{equation}\label{eq:QV-7}
\frac{1}{G_1(x^{-1}R,R)}\ge(2B_t(1))^{-2}G_1(R,x^{-1}R).
\end{equation}
%%%
From (\ref{eq:QV-lower-4}) and (\ref{eq:QV-6}) we obtain for $x\in(0,1]$,
%%%
\begin{eqnarray}
\label{eq:QV-8}
(2B_t(1))^{-2}G_1(xR,R)
&\le&
\frac{1}{G_1(R,xR)}
\le
C_t^2\frac{G_w(xR,R)}{H_w(xR,R)}
\\
&=&
C_t^2G_w(xR,R)H_w(R,xR).
\nonumber
\end{eqnarray}
%%%
Similarly, from (\ref{eq:QV-lower-6}) and (\ref{eq:QV-7}) we obtain for
$x\in(1,\infty)$,
%%%
\begin{equation}\label{eq:QV-9}
(2B_t(1))^{-2}G_1(R,x^{-1}R)\le C_t^2G_w(xR,R)H_w(x^{-1}R,R).
\end{equation}
%%%
Taking the supremum over $R\in(0,d_t]$ in (\ref{eq:QV-8}) and (\ref{eq:QV-9}),
we get
\[
(Q_t1)(x)\le (2C_tB_t(1))^2(Q_tw)(x)(V_tw)(x^{-1}),
\quad
x\in(0,\infty).
\]
From this inequality it follows that for $x\in(0,1]$,
%%%
\begin{equation}\label{eq:QV-10}
\frac{\log (Q_t1)(x)}{\log x}
\ge
\frac{\log (2C_tB_t(1))^2}{\log x}
+
\frac{\log (Q_tw)(x)}{\log x}
-
\frac{\log (V_tw)(x^{-1})}{\log x^{-1}}
\end{equation}
%%%
and, analogously, for $x\in(1,\infty)$,
%%%
\begin{equation}\label{eq:QV-11}
\frac{\log (Q_t1)(x)}{\log x}
\le
\frac{\log (2C_tB_t(1))^2}{\log x}
+
\frac{\log (Q_tw)(x)}{\log x}
-
\frac{\log (V_tw)(x^{-1})}{\log x^{-1}}.
\end{equation}
%%%
Passing to the limit in (\ref{eq:QV-10}) as $x\to 0$ and in (\ref{eq:QV-11})
as $x\to\infty$, we obtain, respectively,
%%%
\begin{equation}\label{eq:QV-12}
\alpha(Q_t1)\ge\alpha(Q_tw)-\beta(V_tw),
\quad
\beta(Q_t1)\le\beta(Q_tw)-\alpha(V_tw).
\end{equation}
%%%

By Lemma~\ref{le:BMO-aux}(a), there exist constants $C_1(t),C_2(t)>0$
such that for every $R>0$,
%%%
\begin{eqnarray}
&&
\exp(-\Omega_t(\log w,R))
\frac{\|w\chi_{\Gamma(t,R)}\|_X\|\chi_{\Gamma(t,R)}\|_{X'}}{|\Gamma(t,R)|}
\le C_1(t),
\label{eq:QV-13}
\\
&&
\exp(\Omega_t(\log w,R))
\frac{\|\chi_{\Gamma(t,R)}\|_X\|\chi_{\Gamma(t,R)}/w\|_{X'}}{|\Gamma(t,R)|}
\le C_2(t).
\label{eq:QV-14}
\end{eqnarray}
%%%
On the other hand, from the lattice property, the H\"older inequality
(see Lemma~\ref{le:Hoelder}), (\ref{eq:loc-Carleson}) and
(\ref{eq:Delta-lower}) it follows that for $R\in(0,d_t]$,
%%%
\begin{eqnarray}
\label{eq:QV-15}
\frac{|\Gamma(t,R)|}{\|\chi_{\Gamma(t,R)}\|_{X'}}
&\le&
\frac{|\Gamma(t,R)|}{\|\chi_{\Delta(t,R)}\|_{X'}}
=
\frac{|\Gamma(t,R)|\cdot\|\chi_{\Delta(t,R)}\|_X}
     {\|\chi_{\Delta(t,R)}\|_X\|\chi_{\Delta(t,R)}\|_{X'}}
\\
&\le&
\frac{|\Gamma(t,R)|}{|\Delta(t,R)|}\|\chi_{\Delta(t,R)}\|_X
\le
\frac{C_{\Gamma,t}R}{R/2}\|\chi_{\Delta(t,R)}\|_X
\nonumber\\
&=&
2C_{\Gamma,t}\|\chi_{\Delta(t,R)}\|_X.
\nonumber
\end{eqnarray}
%%%
Analogously, for $R\in(0,d_t]$,
%%%
\begin{equation}\label{eq:QV-16}
\frac{|\Gamma(t,R)|}{\|\chi_{\Gamma(t,R)}\|_X}
\le
2C_{\Gamma,t}\|\chi_{\Delta(t,R)}\|_{X'}.
\end{equation}
%%%
From (\ref{eq:QV-13})--(\ref{eq:QV-16}) and the lattice property
it follows that for $R\in(0,d_t]$ and $x\in(0,1]$,
%%%
\begin{eqnarray}
\label{eq:QV-17}
&&
G_w(xR,R)
=
\frac{\|w\chi_{\Delta(t,xR)}\|_X\|\chi_{\Delta(t,R)}/w\|_{X'}}{|\Delta(t,R)|}
\\
&&\le
\frac{\|w\chi_{\Gamma(t,xR)}\|_X\|\chi_{\Gamma(t,R)}/w\|_{X'}}{|\Delta(t,R)|}
\nonumber\\
&&\le
\frac{C_1(t)C_2(t)}{|\Delta(t,R)|}
\exp(\Omega_t(\log w,xR))\exp(-\Omega_t(\log w,R))
\nonumber\\
&&\times
\frac{|\Gamma(t,xR)|}{\|\chi_{\Gamma(t,xR)}\|_{X'}}
\cdot
\frac{|\Gamma(t,R)|}{\|\chi_{\Gamma(t,R)}\|_{X}}
\nonumber\\
&&\le
(2C_{\Gamma,t})^2C_1(t)C_2(t) H_w(xR,R)
\frac{\|\chi_{\Delta(t,xR)}\|_X\|\chi_{\Delta(t,R)}\|_{X'}}{|\Delta(t,R)|}
\nonumber\\
&&=
(2C_{\Gamma,t})^2C_1(t)C_2(t) H_w(xR,R) G_1(xR,R)
\nonumber
\end{eqnarray}
%%%
and, similarly, for $R\in(0,d_t]$ and $x\in(1,\infty)$,
%%%
\begin{equation}\label{eq:QV-18}
G_w(R,x^{-1}R)\le
(2C_{\Gamma,t})^2C_1(t)C_2(t) H_w(R,x^{-1}R) G_1(R,x^{-1}R).
\end{equation}
%%%
Taking the supremum over all $R\in(0,d_t]$ in (\ref{eq:QV-17}) and (\ref{eq:QV-18}),
we obtain
\[
(Q_tw)(x)\le
(2C_{\Gamma,t})^2C_1(t)C_2(t)
(V_tw)(x)(Q_t1)(x),
\quad x\in(0,\infty).
\]
Therefore,
%%%
\begin{equation}\label{eq:QV-19}
\alpha(Q_tw)\ge\alpha(V_tw)+\alpha(Q_t1),
\quad
\beta(Q_tw)\le\beta(V_tw)+\beta(Q_t1).
\end{equation}
%%%
Combining (\ref{eq:QV-3}), (\ref{eq:QV-12}), and (\ref{eq:QV-19}),
we arrive at (\ref{eq:QV-1})--(\ref{eq:QV-2}).
\end{proof}
%%%%%%%%%%%%%%%%%%%%%%%%%%%%%%%%%%%%%%%%%%%%%%%%%%%%%%%%%%%%%%%%%%%%%%%%%%%%
If $X$ is a rearrangement-invariant Banach function space, then from
(\ref{eq:RI-fundamental}) it follows that the conditions $1\in A_X(\Gamma,t)$ and
$1\in A_X(\Gamma)$ are equivalent to (\ref{eq:loc-Carleson}) and
(\ref{eq:Carleson}), respectively. Hence, $w\in A_X(\Gamma,t)$
implies $1\in A_X(\Gamma,t)$ whenever $X$ is rearrangement-invariant.
This property allows us to simplify the formulation of Theorem~\ref{th:QV}
for rearrangement-invariant Banach function spaces
(see \cite[Theorems~2.6 and 2.7]{K00}).

Note that $\alpha(Q_t1)$ and $\beta(Q_t1)$ can be considered as a
generalization of the Zippin (fundamental) indices $p_X$ and $q_X$ of a
rearrangement-invariant Banach function space $X$ \cite{Zippin71}.
If $X$ is rearrangement-invariant, then $\alpha(Q_t1)=p_X$
and $\beta(Q_t1)=q_X$ (see \cite[Theorem~5.4]{K98}). On the other hand,
the Zippin indices for an Orlicz space $L^\varphi$ coincide with the reciprocals
of the Matuszewska-Orlicz indices, which control the growth of the Young
function $\varphi$ (see, e.g., \cite{Maligranda85} and the references
given there). The notion of Matuszewska-Orlicz indices of Orlicz spaces
was extended to the case of Musielak-Orlicz spaces in \cite{Kaminska98,KT90}.
Remind that Orlicz spaces are always rearrangement-invariant,
but Musielak-Orlicz spaces are not rearrangement-invariant, in general.
%%%%%%%%%%%%%%%%%%%%%%%%%%%%%%%%%%%%%%%%%%%%%%%%%%%%%%%%%%%%%%%%%%%%%%%%%%%
\subsection{Case of Nakano spaces}
Suppose $\Gamma$ is a rectifiable Jordan curve. Assume that $p:\Gamma\to(1,\infty)$
is a continuous function. Then
%%%
\begin{equation}\label{eq:p-cont-reflexive}
1<p_*:=\min_{t\in\Gamma} p(t)\le\max_{t\in\Gamma} p(t):=p^*<\infty,
\end{equation}
%%%
due to the compactness of $\Gamma$. We will say that a continuous
function $p:\Gamma\to(1,\infty)$ belongs to the class $\cP_t$ if
there is a constant $A_t>0$ such that
%%%
\begin{equation}\label{eq:loc-smooth}
|p(\tau)-p(t)|\le \frac{A_t}{-\log|\tau-t|}
\quad\mbox{for all}\quad \tau\in\Gamma(t,1/2).
\end{equation}
%%%
The class of all continuous functions $p:\Gamma\to(1,\infty)$ such
that $p\in\cP_t$ for every $t\in\Gamma$ and
\[
\sup_{t\in\Gamma}A_t=:A<\infty
\]
is denoted by $\cP$. Clearly, $\cP\subset\cP_t$ for every $t\in\Gamma$.

The class $\cP$ plays a very important role in questions on the boundedness
of maximal functions and singular integrals on (weighted) Nakano spaces
(see \cite{Diening02,KS02-3,PR01}, the references therein,
and also Theorem~\ref{th:Samko-Kokilashvili}).
%%%%%%%%%%%%%%%%%%%%%%%%%%%%%%%%%%%%%%%%%%%%%%%%%%%%%%%%%%%%%%%%%%%%%%%%
\begin{proposition}\label{pr:p-duality}
A function $p$ belongs to $\cP_t$
(respectively, to $\cP$) if and only if the function
$p'(\tau):=p(\tau)/(p(\tau)-1)$ belongs to $\cP_t$ (respectively, to $\cP$).
\end{proposition}
%%%%%%%%%%%%%%%%%%%%%%%%%%%%%%%%%%%%%%%%%%%%%%%%%%%%%%%%%%%%%%%%%%%%%%%
\begin{proof}
The statement immediately follows from the obvious inequality
\[
|p'(\tau)-p'(t)|=\left|\frac{p(\tau)-p(t)}{(p(\tau)-1)(p(t)-1)}\right|
\le
\frac{|p(\tau)-p(t)|}{(p_*-1)^2},
\quad \tau,t\in\Gamma,
\]
and the reflexive relation $(p')'=p$.
\end{proof}
%%%%%%%%%%%%%%%%%%%%%%%%%%%%%%%%%%%%%%%%%%%%%%%%%%%%%%%%%%%%%%%%%%%%%%%%
\begin{lemma}\label{le:ch}
Let $\Gamma$ be locally a Carleson curve at $t\in\Gamma$ and $p\in\cP_t$.
Then there exist constants $M_1(t), M_2(t),C_1(t),C_2(t)\in(0,\infty)$
such that
%%%
\begin{eqnarray}
\|\chi_{\Delta(t,R)}\|_{L^{p(\cdot)}}
&\ge&
M_1(t)R^{1/p(t)}
\quad\mbox{for all}
\quad R\in(0,C_1(t)),
\label{eq:character-1}
\\
\|\chi_{\Gamma(t,R)}\|_{L^{p(\cdot)}}
&\le &
M_2(t)R^{1/p(t)}
\quad\mbox{for all}
\quad R\in(0,C_2(t)).
\label{eq:character-2}
\end{eqnarray}
\end{lemma}
%%%%%%%%%%%%%%%%%%%%%%%%%%%%%%%%%%%%%%%%%%%%%%%%%%%%%%%%%%%%%%%%%%%%%%%%%
\begin{proof}
From (\ref{eq:loc-smooth}) it follows that
%%%
\begin{equation}\label{eq:character-3}
-p(t)-\frac{A_t}{-\log|\tau-t|}\le -p(\tau)\le -p(t)+\frac{A_t}{-\log|\tau-t|},
\quad \tau\in\Gamma(t,1/2).
\end{equation}
%%%
Since $|\tau-t|\le R$ for $\tau\in\Gamma(t,R)$, we have
%%%
\begin{equation}\label{eq:character-4}
\frac{A_t}{-\log|\tau-t|}\le \frac{A_t}{-\log R}, \quad\tau\in\Gamma(t,R), \quad R\in(0,1/2).
\end{equation}
%%%
From (\ref{eq:character-3}) and (\ref{eq:character-4}) we get for
$\tau\in\Gamma(t,R)$ and $R\in(0,1/2)$,
%%%
\begin{equation}\label{eq:character-5}
-p(t)+\frac{A_t}{\log R}\le -p(\tau)\le -p(t)-\frac{A_t}{\log R}.
\end{equation}
%%%
For $R\in(0,e^{-A_t})$, taking into account that $p(t)\in(1,\infty)$, we obtain
%%%
\begin{equation}\label{eq:character-5*}
p(t)+\frac{A_t}{\log R}=(p(t)-1)+\left(1+\frac{A_t}{\log R}\right)>p(t)-1>0.
\end{equation}
%%%
From (\ref{eq:character-5}) we get for $\lambda\in(0,1]$ and $R\in(0,\min\{1/2,e^{-A_t}\})$,
%%%
\begin{eqnarray}
\label{eq:character-6}
\exp\left(-\left[p(t)+\frac{A_t}{\log R}\right]\log\lambda\right)
&\le&
\exp(-p(\tau)\log\lambda)
\\
&\le&
\exp\left(-\left[p(t)-\frac{A_t}{\log R}\right]\log\lambda\right).
\nonumber
\end{eqnarray}
%%%
Analogously, for $\lambda\in(1,\infty)$ and $R\in(0,\min\{1/2,e^{-A_t}\})$,
%%%
\begin{eqnarray}
\label{eq:character-7}
\exp\left(-\left[p(t)-\frac{A_t}{\log R}\right]\log\lambda\right)
&\le&
\exp(-p(\tau)\log\lambda)
\\
&\le&
\exp\left(-\left[p(t)+\frac{A_t}{\log R}\right]\log\lambda\right).
\nonumber
\end{eqnarray}
%%%

Let us prove (\ref{eq:character-1}). From the first inequality in
(\ref{eq:character-6}) and (\ref{eq:Delta-lower})
it follows that for $\lambda\in(0,1]$ and
$R\in(0,\min\{1/2,e^{-A_t},d_t\})$,
%%%
\begin{eqnarray*}
m(\chi_{\Delta(t,R)}/\lambda,p)
&=&
\int_{\Delta(t,R)}\exp(-p(\tau)\log\lambda)|d\tau|
\\
&\ge&
\exp\left(-\left[p(t)+\frac{A_t}{\log R}\right]\log\lambda\right)|\Delta(t,R)|
\\
&\ge&
\exp\left(\log\frac{R}{2}-\left[p(t)+\frac{A_t}{\log R}\right]\log\lambda\right).
\end{eqnarray*}
%%%
Put $C_1(t):=\min\{1/2, e^{-A_t},d_t\}$.
Therefore, taking into account (\ref{eq:character-5*}), we obtain for
$R\in(0,C_1(t))$,
%%%
\begin{eqnarray*}
&&
\Big\{\lambda\in(0,1]: \quad m(\chi_{\Delta(t,R)}/\lambda,p)\le 1\Big\}
\\
&&
\subset
\left\{\lambda\in(0,1]: \quad
\log\frac{R}{2}-\left[p(t)+\frac{A_t}{\log R}\right]\log\lambda\le 0
\right\}
\\
&&=
\left\{\lambda: \quad
\exp\left(\frac{\log(R/2)}{p(t)+A_t/\log R}\right)\le\lambda\le 1
\right\}.
\end{eqnarray*}
%%%
Thus, for $R\in(0,C_1(t))$,
%%%
\begin{eqnarray}\label{eq:character-8}
&&
N_1 :=
\inf\Big\{\lambda\in(0,1]:\quad m(\chi_{\Delta(t,R)}/\lambda,p)\le 1\Big\}
\ge
\exp\left(\frac{\log(R/2)}{p(t)+A_t/\log R}\right).
\end{eqnarray}
%%%
Analogously, from the first inequality in (\ref{eq:character-7}) we obtain
\[
\Big\{\lambda\in(1,\infty): \quad m(\chi_{\Delta(t,R)}/\lambda,p)\le 1\Big\}
\subset(1,\infty)
\]
because
%%%
\begin{equation}\label{eq:character-9}
\exp\left(\frac{\log(R/2)}{p(t)-A_t/\log R}\right)<1
\quad\mbox{for}\quad R\in(0,C_1(t)).
\end{equation}
%%%
Thus, for $R\in(0,C_1(t))$,
%%%
\begin{equation}\label{eq:character-10}
N_2:=\inf\Big\{\lambda\in(1,\infty): \quad m(\chi_{\Delta(t,R)}/\lambda,p)\le 1\Big\}\ge 1.
\end{equation}
%%%
From (\ref{eq:character-8})--(\ref{eq:character-10}) we obtain for
$R\in(0,C_1(t))$,
%%%
\begin{eqnarray}
\label{eq:character-11}
&&
\|\chi_{\Delta(t,R)}\|_{L^{p(\cdot)}}
=
\inf\Big\{\lambda>0: \quad m(\chi_{\Delta(t,R)}/\lambda,p)\le 1\Big\}
=
\min\{N_1,N_2\}
\\
&&\ge
\min\left\{
1,
\exp\left(\frac{\log(R/2)}{p(t)+A_t/\log R}\right)
\right\}
=
\exp\left(\frac{\log(R/2)}{p(t)+A_t/\log R}\right).
\nonumber
\end{eqnarray}
%%%
From (\ref{eq:character-5*}) it follows that for $R\in(0,C_1(t))$,
\[
\frac{\log\frac{R}{2}}{p(t)+\frac{A_t}{\log R}}
-
\frac{\log\frac{R}{2}}{p(t)}
=
\frac{-A_t+A_t\frac{\log 2}{\log R}}{\left(p(t)+\frac{A_t}{\log R}\right)p(t)}
\ge
\frac{-A_t+A_t\frac{\log 2}{\log R}}{(p(t)-1)p(t)}
\ge
\frac{A_t+\log 2}{(1-p(t))p(t)}.
\]
From the latter inequality we deduce that
%%%
\begin{eqnarray}
\label{eq:character-12}
&&
\exp\left(\frac{\log(R/2)}{p(t)+A_t/\log R}\right)
=
\exp\left(\frac{\log(R/2)}{p(t)+A_t/\log R}-\frac{\log(R/2)}{p(t)}\right)
\left(\frac{R}{2}\right)^{1/p(t)}
\\
&&\ge
\exp\left(\frac{A_t+\log 2}{(1-p(t))p(t)}-\frac{\log 2}{p(t)}\right)R^{1/p(t)}.
\nonumber
\end{eqnarray}
%%%
Combining (\ref{eq:character-11}) and (\ref{eq:character-12}),
we arrive at (\ref{eq:character-1}) with
\[
C_1(t):=\min\{1/2,e^{-A_t},d_t\},
\quad
M_1(t):=\exp\left(\frac{A_t+\log 2}{(1-p(t))p(t)}-\frac{\log 2}{p(t)}\right).
\]
Taking into account (\ref{eq:loc-Carleson}), one can prove
that (\ref{eq:character-2}) is valid with
\[
C_2(t):=\min\{1/2,1/C_{\Gamma,t},e^{-A_t},d_t\},
\quad
M_2(t):=\exp\left(
\frac{A_t}{(p(t))^2}+\frac{\log C_{\Gamma,t}}{p(t)}\right).
\]
The proof of (\ref{eq:character-2}) is similar to the proof of
(\ref{eq:character-1}) and it is omitted.
\end{proof}
%%%%%%%%%%%%%%%%%%%%%%%%%%%%%%%%%%%%%%%%%%%%%%%%%%%%%%%%%%%%%%%%%%%%%%%%
\begin{lemma}\label{le:ALP}
Suppose $\Gamma$ is locally a Carleson curve at $t\in\Gamma$ and $p\in\cP_t$.
Then $1\in A_{L^{p(\cdot)}}(\Gamma,t)$ and
%%%
\begin{equation}\label{eq:ALP-1}
\alpha(Q_t1)=\beta(Q_t1)=1/p(t).
\end{equation}
\end{lemma}
%%%%%%%%%%%%%%%%%%%%%%%%%%%%%%%%%%%%%%%%%%%%%%%%%%%%%%%%%%%%%%%%%%%%%%%%
\begin{proof}
From Lemma~\ref{le:ch} we deduce that there exist constants
$C_i(t), M_i(t)\ (i=1,2)$ such that
%%%
\begin{eqnarray}
\|\chi_{\Delta(t,R)}\|_{L^{p(\cdot)}}
&\ge&
M_1(t)R^{1/p(t)}
\quad\mbox{for all}
\quad R\in(0,C_1(t)),
\label{eq:ALP-2}
\\
\|\chi_{\Gamma(t,R)}\|_{L^{p(\cdot)}}
&\le &
M_2(t)R^{1/p(t)}
\quad\mbox{for all}
\quad R\in(0,C_2(t)),
\label{eq:ALP-3}
\end{eqnarray}
%%%
By Proposition~\ref{pr:p-duality}, $p'\in\cP_t$.
Analogously, applying Lemma~\ref{le:ch} to $L^{p'(\cdot)}$
and taking into account that the latter space coincide with
$(L^{p(\cdot)})'$ up to the equivalence of the norms
(see Lemma~\ref{le:Nakano-duality}), we infer that there
exist constants $C_i'(t), M_i'(t)\ (i=1,2)$ such that
%%%
\begin{eqnarray}
\|\chi_{\Delta(t,R)}\|_{(L^{p(\cdot)})'}
&\ge&
M_1'(t)R^{1/p'(t)}
\quad\mbox{for all}
\quad R\in(0,C_1'(t)),
\label{eq:ALP-4}
\\
\|\chi_{\Gamma(t,R)}\|_{(L^{p(\cdot)})'}
&\le &
M_2'(t)R^{1/p'(t)}
\quad\mbox{for all}
\quad R\in(0,C_2'(t)).
\label{eq:ALP-5}
\end{eqnarray}
%%%
From (\ref{eq:ALP-3}), (\ref{eq:ALP-5}) it follows that
for $R\in(0,\min\{C_2(t),C_2'(t)\})$,
%%%
\begin{eqnarray}
\label{eq:ALP-6}
B_{t,R}(1) &=&
\frac{1}{R}
\|\chi_{\Gamma(t,R)}\|_{L^{p(\cdot)}}
\|\chi_{\Gamma(t,R)}\|_{(L^{p(\cdot)})'}
\\
&\le&
\frac{1}{R}M_2(t)M_2'(t)R^{1/p(t)}R^{1/p'(t)}=M_2(t)M_2'(t).
\nonumber
\end{eqnarray}
%%%
On the other hand, for $R\ge\min\{C_2(t),C_2'(t)\}$,
%%%
\begin{equation}\label{eq:ALP-7}
B_{t,R}(1) =
\frac{1}{R}
\|\chi_{\Gamma(t,R)}\|_{L^{p(\cdot)}}
\|\chi_{\Gamma(t,R)}\|_{(L^{p(\cdot)})'}
\le
\frac{\|1\|_{L^{p(\cdot)}}\|1\|_{(L^{p(\cdot)})'}}{\min\{C_2(t),C_2'(t)\}}.
\end{equation}
%%%
From (\ref{eq:ALP-6}) and (\ref{eq:ALP-7}) it follows that
\[
\sup_{R>0}B_{t,R}(1)\le\max\left\{
M_2(t)M_2'(t),
\frac{\|1\|_{L^{p(\cdot)}}\|1\|_{(L^{p(\cdot)})'}}{\min\{C_2(t),C_2'(t)\}}
\right\}<\infty.
\]
Thus, $1\in A_{L^{p(\cdot)}}(\Gamma,t)$.

Put $C(t):=\min\{C_1(t),C_2(t),C_1'(t),C_2'(t)\}$.
From (\ref{eq:ALP-3}), (\ref{eq:ALP-5}), (\ref{eq:Delta-lower}),
and the lattice property we obtain
for $x\in(0,\infty)$ and $R\in(0,C(t)\min\{1,1/x\})$,
%%%
\begin{eqnarray}
\label{eq:ALP-8}
G_1(xR,R) &:=&
\frac{
\|\chi_{\Delta(t,xR)}\|_{L^{p(\cdot)}}
\|\chi_{\Delta(t, R)}\|_{(L^{p(\cdot)})'}
}{|\Delta(t,R)|}
\\
&\le&
M_2(t)M_2'(t) \frac{(xR)^{1/p(t)}R^{1/p'(t)}}{|\Delta(t,R)|}
\nonumber\\
&\le&
M_2(t)M_2'(t) \frac{x^{1/p(t)}R}{R/2}
=
2M_2(t)M_2'(t)x^{1/p(t)}.
\nonumber
\end{eqnarray}
%%%
Combining (\ref{eq:ALP-4}), (\ref{eq:ALP-6}), and (\ref{eq:Delta-upper}), we get
for the same $x$ and $R$,
%%%
\begin{eqnarray}
\label{eq:ALP-9}
G_1(xR,R)
&\ge&
M_1(t)M_1'(t)\frac{(xR)^{1/p(t)}R^{1/p'(t)}}{|\Delta(t,R)|}
\\
&\ge&
M_1(t)M_1'(t)\frac{x^{1/p(t)}R}{C_{\Gamma,t}R}
=
\frac{M_1(t)M_1'(t)}{C_{\Gamma,t}}x^{1/p(t)}.
\nonumber
\end{eqnarray}
%%%
From (\ref{eq:ALP-8}) and (\ref{eq:ALP-9}) it follows that
\[
\frac{M_1(t)M_1'(t)}{C_{\Gamma,t}}x^{1/p(t)}
\le
(Q_t^01)(x)
\le
2M_2(t)M_2'(t)x^{1/p(t)},
\quad x\in(0,\infty).
\]
Since $1\in A_{L^{p(\cdot)}}(\Gamma)$, the function $Q_t^01$ is
regular and submultiplicative (see Lemmas~\ref{le:Q-sub}
and \ref{le:Q-reg}).
From the latter inequality it follows that
\[
\alpha(Q_t^01)=\beta(Q_t^01)=1/p(t).
\]
Combining the latter equalities with Lemma~\ref{le:Q-sub},
we arrive at (\ref{eq:ALP-1}).
\end{proof}
%%%%%%%%%%%%%%%%%%%%%%%%%%%%%%%%%%%%%%%%%%%%%%%%%%%%%%%%%%%%%%%%%%%%%%%%
\begin{theorem}\label{th:disintegration}
Let $\Gamma$ be locally a Carleson curve at $t\in\Gamma$, let $w:\Gamma\to[0,\infty]$
be a weight, and let $p\in\cP_t$. If $w\in A_{L^{p(\cdot)}}(\Gamma,t)$,
then $\log w\in BMO(\Gamma,t)$ and
\begin{equation}\label{eq:disintegration-1}
\alpha(Q_tw)=1/p(t)+\alpha(V_tw),
\quad
\beta(Q_tw)=1/p(t)+\beta(V_tw).
\end{equation}
\end{theorem}
%%%%%%%%%%%%%%%%%%%%%%%%%%%%%%%%%%%%%%%%%%%%%%%%%%%%%%%%%%%%%%%%%%%%%%%%
\begin{proof}
Since $p\in\cP_t$ and $\Gamma$ is locally a Carleson curve at $t$,
in view of Lemma~\ref{le:ALP}, $1\in A_{L^{p(\cdot)}}(\Gamma,t)$.
By Lemma~\ref{le:BMO}(a), $\log w\in BMO(\Gamma,t)$.
From Theorem~\ref{th:QV} and (\ref{eq:ALP-1}) we get
%%%
\begin{eqnarray*}
1/p(t)+\alpha(V_tw)
\le
\alpha(Q_tw) &\le& \min\{1/p(t)+\alpha(V_tw),1/p(t)+\beta(V_tw)\}
\\
&=&
1/p(t)+\alpha(V_tw),
\\
1/p(t)+\beta(V_tw) \ge\beta(Q_tw) &\ge&
\max\{1/p(t)+\alpha(V_tw),1/p(t)+\beta(V_tw)\}
\\
&=&
1/p(t)+\beta(V_tw),
\end{eqnarray*}
%%%
that is, equalities (\ref{eq:disintegration-1}) hold.
\end{proof}
%%%%%%%%%%%%%%%%%%%%%%%%%%%%%%%%%%%%%%%%%%%%%%%%%%%%%%%%%%%%%%%%%%%%%%%%
\begin{lemma}
Let $\Gamma$ be a Carleson curve, let $w:\Gamma\to[0,\infty]$ be a weight,
and let $p\in\cP$. If $w\in A_{L^{p(\cdot)}}(\Gamma)$, then
$\log w\in BMO(\Gamma)$.
\end{lemma}
%%%%%%%%%%%%%%%%%%%%%%%%%%%%%%%%%%%%%%%%%%%%%%%%%%%%%%%%%%%%%%%%%%%%%%%%
\begin{proof}
By analogy with Lemma~\ref{le:ch} one can show that
there exist constants $C>0$ and $M,M'\in(0,\infty)$ such that
\[
\|\chi_{\Gamma(t,R)}\|_{L^{p(\cdot)}}\le MR^{1/p(t)},
\quad
\|\chi_{\Gamma(t,R)}\|_{L^{p'(\cdot)}}\le M'R^{1/p'(t)}
\]
for all $R\in(0,C)$ and all $t\in\Gamma$. Taking into account
Lemma~\ref{le:Nakano-duality}, as in Lemma~\ref{le:ALP}
from the latter inequalities we obtain $1\in A_{L^{p(\cdot)}}(\Gamma)$.
Therefore, $\log w\in BMO(\Gamma)$, due to Lemma~\ref{le:BMO}(b).
\end{proof}
%%%%%%%%%%%%%%%%%%%%%%%%%%%%%%%%%%%%%%%%%%%%%%%%%%%%%%%%%%%%%%%%%%%%%%%%
\subsection{Indicator functions}
In this subsection we generalize the notion of indicator functions
(see \cite[Ch.~3]{BK97} and also \cite[Section~7.2]{K98}, \cite[Section~2.5]{K00},
\cite[Section~3.3]{K02})
to the case of weighted Banach function spaces.

Suppose $\Gamma$ is a rectifiable Jordan curve, $w:\Gamma\to[0,\infty]$
is a weight, $X$ is a Banach function space.
%%%%%%%%%%%%%%%%%%%%%%%%%%%%%%%%%%%%%%%%%%%%%%%%%%%%%%%%%%%%%%%%%%%%%%%%
\begin{lemma}\label{le:indicator0}
Let $\Gamma$ be locally a Carleson curve at $t\in\Gamma$. For every
$x\in\R$, the function $W_t\eta_t^x$ is regular, submultiplicative, and
\begin{eqnarray*}
\alpha_t^0(x) &:=&\alpha(W_t\eta_t^x)=\min\{\delta_t^-x,\delta_t^+x\},\\
\beta_t^0(x) &:=&\beta(W_t\eta_t^x)=\max\{\delta_t^-x,\delta_t^+x\}.
\end{eqnarray*}
\end{lemma}

This statement follows from local analogs of \cite[Lemmas~1.15, 1.16, and
Proposition~3.1]{BK97}.

For a complex number $\gamma\in\C$, we define a continuous
function $\varphi_{t,\gamma}$ on $\Gamma\setminus\{t\}$ by
%%%
\begin{equation}\label{eq:canonical-weight}
\varphi_{t,\gamma}(\tau)
:=
|(\tau-t)^\gamma|
=
|\tau-t|^{\operatorname{Re}\gamma}
e^{-\operatorname{Im}\gamma\arg(\tau-t)}
=
|\tau-t|^{\operatorname{Re}\gamma}
(\eta_t(\tau))^{\operatorname{Im}\gamma}.
\end{equation}
%%%%%%%%%%%%%%%%%%%%%%%%%%%%%%%%%%%%%%%%%%%%%%%%%%%%%%%%%%%%%%%%%%%%%%%%
\begin{lemma}\label{le:indicator*}
If $w\in A_X(\Gamma,t)$, then for every $\gamma\in\C$, the function
$Q_t(\varphi_{t,\gamma}w)$ is regular, submultiplicative, and
%%%
\begin{eqnarray}
\alpha(Q_t(\varphi_{t,\gamma}w))
&=&
\operatorname{Re}\gamma+\alpha(Q_t(\eta_t^{\operatorname{Im}\gamma}w)),
\label{eq:indicator*-1}
\\
\beta(Q_t(\varphi_{t,\gamma}w))
&=&
\operatorname{Re}\gamma+\beta(Q_t(\eta_t^{\operatorname{Im}\gamma}w)).
\label{eq:indicator*-2}
\end{eqnarray}
\end{lemma}
%%%%%%%%%%%%%%%%%%%%%%%%%%%%%%%%%%%%%%%%%%%%%%%%%%%%%%%%%%%%%%%%%%%%%%%%
\begin{proof}
This statement is proved similarly to \cite[Lemma~7.2]{K98}.
By a local analog of \cite[Proposition~3.1]{BK97}, the function
$W_t\varphi_{t,\operatorname{Re}\gamma}$
is regular and submultiplicative for every $\gamma\in\C$ and
%%%
\begin{equation}\label{eq:indicator*-3}
\alpha(W_t\varphi_{t,\operatorname{Re}\gamma})
=
\beta(W_t\varphi_{t,\operatorname{Re}\gamma})
=
\operatorname{Re}\gamma.
\end{equation}
%%%
On the other hand, by Lemmas~\ref{le:Q-sub}--\ref{le:Q-reg},
the function $Q_tw$ is regular and submultiplicative. Then,
by Theorem~\ref{th:QW}, the function $Q_t(\varphi_{t,\gamma}w)$
is regular and submultiplicative for every $\gamma\in\C$.
In particular, the function $Q_t(\eta_t^{\operatorname{Im}\gamma}w)$
is regular and submultiplicative for every $\gamma\in\C$.
From Theorem~\ref{th:QW} and (\ref{eq:indicator*-3}) it follows
that
%%%
\begin{eqnarray*}
\alpha(Q_t(\eta_t^{\operatorname{Im}\gamma}w))+\operatorname{Re}\gamma
&\le&
\alpha(Q_t(\varphi_{t,\gamma}w))
\\
&\le&
\min\{
\alpha(Q_t(\eta_t^{\operatorname{Im}\gamma}w))+\operatorname{Re}\gamma,
\beta(Q_t(\eta_t^{\operatorname{Im}\gamma}w))+\operatorname{Re}\gamma
\},
\\
\beta(Q_t(\eta_t^{\operatorname{Im}\gamma}w))+\operatorname{Re}\gamma
&\ge&
\beta(Q_t(\varphi_{t,\gamma}w))
\\
&\ge&
\max\{
\alpha(Q_t(\eta_t^{\operatorname{Im}\gamma}w))+\operatorname{Re}\gamma,
\beta(Q_t(\eta_t^{\operatorname{Im}\gamma}w))+\operatorname{Re}\gamma
\}.
\end{eqnarray*}
%%%
From the latter inequalities we immediately obtain
(\ref{eq:indicator*-1})--(\ref{eq:indicator*-2}).
\end{proof}
%%%%%%%%%%%%%%%%%%%%%%%%%%%%%%%%%%%%%%%%%%%%%%%%%%%%%%%%%%%%%%%%%%%%%%%%
\begin{lemma}\label{le:indicator}
If $w\in A_X(\Gamma,t)$ and $1\in A_X(\Gamma,t)$, then for every $\gamma\in\C$,
the function $V_t(\varphi_{t,\gamma}w)$ is regular, submultiplicative, and
%%%
\begin{eqnarray}
\alpha(V_t(\varphi_{t,\gamma}w))
&=&
\operatorname{Re}\gamma+\alpha(V_t(\eta_t^{\operatorname{Im}\gamma}w)),
\label{eq:indicator-1}
\\
\beta(V_t(\varphi_{t,\gamma}w))
&=&
\operatorname{Re}\gamma+\beta(V_t(\eta_t^{\operatorname{Im}\gamma}w)).
\label{eq:indicator-2}
\end{eqnarray}
\end{lemma}
%%%%%%%%%%%%%%%%%%%%%%%%%%%%%%%%%%%%%%%%%%%%%%%%%%%%%%%%%%%%%%%%%%%%%%%%
\begin{proof}
By Lemma~\ref{le:BMO}(a), $\log w\in BMO(\Gamma,t)$. Then by
Lemma~\ref{le:V-reg}, the function $V_tw$ is regular. The rest is proved
by analogy with Lemma~\ref{le:indicator*} with the help of
Theorem~\ref{th:VW}.
\end{proof}
%%%%%%%%%%%%%%%%%%%%%%%%%%%%%%%%%%%%%%%%%%%%%%%%%%%%%%%%%%%%%%%%%%%%%%%%
If $w\in A_X(\Gamma,t)$, then for every $x\in\R$, the function
$Q_t(\eta_t^xw)$ is regular and submultiplicative, in view of
Lemma~\ref{le:indicator*}. From Theorem~\ref{th:submult} and
Lemma~\ref{le:Q-reg} we deduce that the following functions
are well defined for $x\in\R$:
\[
\alpha_t^*(x):=\alpha(Q_t(\eta_t^xw))=\alpha(Q_t^0(\eta_t^xw)),
\quad
\beta_t^*(x):=\beta(Q_t(\eta_t^xw))=\beta(Q_t^0(\eta_t^xw)).
\]
If, in addition, $1\in A_X(\Gamma,t)$, then the function $V_t(\eta_t^xw)$
is regular and submultiplicative for each $x\in\R$, due to
Lemma~\ref{le:indicator}. Then Theorem~\ref{th:submult}
and Lemma~\ref{le:V-sub} imply that the functions
\[
\alpha_t(x):=\alpha(V_t(\eta_t^xw))=\alpha(V_t^0(\eta_t^xw)),
\quad
\beta_t(x):=\beta(V_t(\eta_t^xw))=\beta(V_t^0(\eta_t^xw))
\]
are well defined for all $x\in\R$.

The functions $\alpha_t^*,\beta_t^*$ are called the
\textit{indicator functions of the triple} $(\Gamma,X,w)$ \textit{at}
$t\in\Gamma$. The functions $\alpha_t,\beta_t$ are referred to as the
\textit{indicator functions of the pair} $(\Gamma,w)$ \textit{at}
$t\in\Gamma$. The functions $\alpha_t^*,\beta_t^*$
were introduced in \cite{K00} (see also \cite{K98,K02})
for rearrangement-invariant Banach function spaces.
The functions $\alpha_t,\beta_t$ were defined
in \cite[Ch.~3]{BK97} in the context of Lebesgue spaces and
Muckenhoupt weights.
%%%%%%%%%%%%%%%%%%%%%%%%%%%%%%%%%%%%%%%%%%%%%%%%%%%%%%%%%%%%%%%%%%%%%%%%
\begin{lemma}
The functions $\alpha_t,\alpha_t^*$ are concave, the functions
$\beta_t,\beta_t^*$ are convex. In particular, $\alpha_t,\alpha_t^*$
and $\beta_t,\beta_t^*$ are continuous on $\R$.
\end{lemma}
%%%%%%%%%%%%%%%%%%%%%%%%%%%%%%%%%%%%%%%%%%%%%%%%%%%%%%%%%%%%%%%%%%%%%%%%
\begin{proof}
By \cite[Section~2.2, Property 6]{KPS78},
\[
\Big\| |f|^\theta|g|^{1-\theta}\Big\|_X
\le
\|f\|_X^\theta\|g\|_X^{1-\theta},
\quad\theta\in[0,1],
\]
for every $f,g\in X$. With the help of this property, one can prove concavity
of $\alpha_t^*$ and convexity of $\beta_t^*$ similarly to
\cite[Proposition~3.20]{BK97}. Concavity of $\alpha_t$
and convexity of $\beta_t$ are already proved there.
\end{proof}
%%%%%%%%%%%%%%%%%%%%%%%%%%%%%%%%%%%%%%%%%%%%%%%%%%%%%%%%%%%%%%%%%%%%%%%%
The following statement generalizes \cite[Lemma~3.5]{K02}.
%%%%%%%%%%%%%%%%%%%%%%%%%%%%%%%%%%%%%%%%%%%%%%%%%%%%%%%%%%%%%%%%%%%%%%%%
\begin{lemma}\label{le:ind-relations}
{\rm (a)} If $w\in A_X(\Gamma,t)$, then for $x,y\in\R$,
\begin{eqnarray*}
\alpha_t^*(x)+\alpha_t^0(y)
\le
&\alpha_t^*(x+y)&
\le
\min\{\alpha_t^*(x)+\beta_t^0(y),\beta_t^*(x)+\alpha_t^0(y)\},
\\
\beta_t^*(x)+\beta_t^0(y)
\ge
&\beta_t^*(x+y)&
\ge
\max\{\alpha_t^*(x)+\beta_t^0(y),\beta_t^*(x)+\alpha_t^0(y)\}.
\end{eqnarray*}

\noindent
{\rm (b)} If $w\in A_X(\Gamma,t)$ and $1\in A_X(\Gamma,t)$, then for $x,y\in\R$,
\begin{eqnarray*}
\alpha_t(x)+\alpha_t^0(y)
\le
&\alpha_t(x+y)&
\le
\min\{\alpha_t(x)+\beta_t^0(y),\beta_t(x)+\alpha_t^0(y)\},
\\
\beta_t(x)+\beta_t^0(y)
\ge
&\beta_t(x+y)&
\ge
\max\{\alpha_t(x)+\beta_t^0(y),\beta_t(x)+\alpha_t^0(y)\}.
\end{eqnarray*}
\end{lemma}
%%%%%%%%%%%%%%%%%%%%%%%%%%%%%%%%%%%%%%%%%%%%%%%%%%%%%%%%%%%%%%%%%%%%%%%%
\begin{proof}
(a) From Lemmas~\ref{le:indicator0} and~\ref{le:indicator*}
it follows that the functions $\alpha_t^*,\beta_t^*$ and
$\alpha_t^0,\beta_t^0$ are well defined.
Applying Theorem~\ref{th:QW} to the weights $w:=\eta_t^xw$ and
$\psi:=\eta_t^y$, we get Part (a).
Part (b) is proved analogously with the help of Theorem~\ref{th:VW}
and Lemma~\ref{le:indicator}.
\end{proof}
%%%%%%%%%%%%%%%%%%%%%%%%%%%%%%%%%%%%%%%%%%%%%%%%%%%%%%%%%%%%%%%%%%%%%%%%
\begin{corollary}\label{co:ind-relations}
Let $\Gamma$ be locally a Carleson curve at $t\in\Gamma$ such that
$\delta_t^-=\delta_t^+=:\delta_t$.

\noindent
{\rm (a)} If $w\in A_X(\Gamma,t)$, then
\begin{equation}\label{eq:ind-relations-1}
\alpha_t^*(x)=\alpha(Q_tw)+\delta_t x,
\quad
\beta_t^*(x)=\beta(Q_tw)+\delta_t x
\quad (x\in\R).
\end{equation}

\noindent
{\rm (b)} If $w\in A_X(\Gamma,t)$ and $1\in A_X(\Gamma,t)$, then
\begin{equation}\label{eq:ind-relations-2}
\alpha_t(x)=\mu_t+\delta_t x,
\quad
\beta_t(x)=\nu_t+\delta_t x
\quad (x\in\R).
\end{equation}
\end{corollary}
%%%%%%%%%%%%%%%%%%%%%%%%%%%%%%%%%%%%%%%%%%%%%%%%%%%%%%%%%%%%%%%%%%%%%%%%
\begin{proof}
(a) Since $\delta_t^-=\delta_t^+=\delta_t$, we have
$\alpha_t^0(x)=\beta_t^0(x)=\delta_tx$. In that case from
Lemma~\ref{le:ind-relations}(a) we deduce that
%%%
\begin{equation}\label{eq:ind-relations-3}
\alpha_t^*(y)+\delta_tx=\alpha_t^*(x+y),
\quad
\beta_t^*(y)+\delta_tx=\beta_t^*(x+y)
\end{equation}
%%%
for every $x,y\in\R$. Setting $y=0$ in (\ref{eq:ind-relations-3}),
we arrive at (\ref{eq:ind-relations-1}). Part (b) is proved similarly.
\end{proof}
%%%%%%%%%%%%%%%%%%%%%%%%%%%%%%%%%%%%%%%%%%%%%%%%%%%%%%%%%%%%%%%%%%%%%%%%
\subsection{Indicator functions for Nakano spaces}
Let $\Gamma$ be a rectifiable Jordan curve, let $L^{p(\cdot)}$ be a
Nakano space. Fix $t\in\Gamma$. For a weight
$w\in A_{L^{p(\cdot)}}(\Gamma,t)$, put
\[
N_t:=\Big\{
\gamma\in\C: \quad \varphi_{t,\gamma}w\in A_{L^{p(\cdot)}}(\Gamma,t)
\Big\}.
\]
%%%%%%%%%%%%%%%%%%%%%%%%%%%%%%%%%%%%%%%%%%%%%%%%%%%%%%%%%%%%%%%%%%%%%%%%
\begin{lemma}\label{le:Nakano-ind-Carleson}
Let $\Gamma$ be locally a Carleson curve at $t\in\Gamma$, let $p\in\cP_t$,
and let $w\in A_{L^{p(\cdot)}}(\Gamma,t)$.
Then for every $\gamma\in N_t$,
%%%
\begin{equation}\label{eq:Nakano-ind-Carleson-1}
\alpha_t^*(\operatorname{Im}\gamma)=1/p(t)+\alpha_t(\operatorname{Im}\gamma),
\quad
\beta_t^*(\operatorname{Im}\gamma)=1/p(t)+\beta_t(\operatorname{Im}\gamma).
\end{equation}
\end{lemma}
%%%%%%%%%%%%%%%%%%%%%%%%%%%%%%%%%%%%%%%%%%%%%%%%%%%%%%%%%%%%%%%%%%%%%%%%
\begin{proof}
Let $\gamma\in N_t$. By Theorem~\ref{th:disintegration},
%%%
\begin{eqnarray}
\alpha(Q_t(\varphi_{t,\gamma}w))
&=&
1/p(t)+\alpha(V_t(\varphi_{t,\gamma}w)),
\label{eq:Nakano-ind-Carleson-2}
\\
\beta(Q_t(\varphi_{t,\gamma}w))
&=&
1/p(t)+\beta(V_t(\varphi_{t,\gamma}w)).
\label{eq:Nakano-ind-Carleson-3}
\end{eqnarray}
%%%
Note that by Lemma~\ref{le:ALP}, $1\in A_{L^{p(\cdot)}}(\Gamma,t)$.
Therefore, we can apply Lemma~\ref{le:indicator}. From
(\ref{eq:Nakano-ind-Carleson-2})--(\ref{eq:Nakano-ind-Carleson-3}),
(\ref{eq:indicator*-1})--(\ref{eq:indicator*-2}), and
(\ref{eq:indicator-1})--(\ref{eq:indicator-2}) it follows that
\[
\alpha(Q_t(\eta_t^{\operatorname{Im}\gamma}w))
=
1/p(t)+\alpha(V_t(\eta_t^{\operatorname{Im}\gamma}w)),
\quad
\beta(Q_t(\eta_t^{\operatorname{Im}\gamma}w))
=
1/p(t)+\beta(V_t(\eta_t^{\operatorname{Im}\gamma}w)),
\]
that is, equalities (\ref{eq:Nakano-ind-Carleson-1}) hold.
\end{proof}
%%%%%%%%%%%%%%%%%%%%%%%%%%%%%%%%%%%%%%%%%%%%%%%%%%%%%%%%%%%%%%%%%%%%%%%%
\begin{lemma}\label{le:Nakano-ind-smooth}
Let $\Gamma$ be locally a Carleson  curve at $t\in\Gamma$
such that $\delta_t^-=\delta_t^+=0$, let $p\in\cP_t$, and
let  $w\in A_{L^{p(\cdot)}}(\Gamma,t)$. Then for every  every $x\in\R$,
\begin{equation}\label{eq:Nakano-ind-smooth-1}
\alpha_t(x)=\mu_t,
\quad
\beta_t(x)=\nu_t,
\quad
\alpha_t^*(x)=1/p(t)+\mu_t,
\quad
\beta_t^*(x)=1/p(t)+\nu_t,
\end{equation}
%%%
where $\mu_t,\nu_t$ are the indices of powerlikeness of the
weight $w$ at $t$ defined by {\rm (\ref{eq:ind_powerlikeness})}.
\end{lemma}
%%%%%%%%%%%%%%%%%%%%%%%%%%%%%%%%%%%%%%%%%%%%%%%%%%%%%%%%%%%%%%%%%%%%%%%%
\begin{proof}
By Lemma~\ref{le:ALP}, $1\in A_{L^{p(\cdot)}}(\Gamma,t)$.
From Corollary~\ref{co:ind-relations} we get for every $x\in\R$,
%%%
\begin{equation}\label{eq:Nakano-ind-smooth-2}
\alpha_t^*(x)=\alpha(Q_tw),
\quad
\beta_t^*(x)=\beta_t(Q_tw),
\quad
\alpha_t(x)=\mu_t,
\quad
\beta_t(x)=\nu_t.
\end{equation}
%%%
On the other hand, by Theorem~\ref{th:disintegration},
%%%
\begin{equation}\label{eq:Nakano-ind-smooth-3}
\alpha(Q_tw)=1/p(t)+\mu_t,\quad\beta(Q_tw)=1/p(t)+\nu_t.
\end{equation}
%%%
Combining (\ref{eq:Nakano-ind-smooth-2}) and (\ref{eq:Nakano-ind-smooth-3}),
we arrive at (\ref{eq:Nakano-ind-smooth-1}).
\end{proof}
%%%%%%%%%%%%%%%%%%%%%%%%%%%%%%%%%%%%%%%%%%%%%%%%%%%%%%%%%%%%%%%%%%%%%%%%
\section{Fredholm theory for singular integral operators\\
with bounded measurable coefficients}\label{sec:bounded}
\subsection{The Cauchy singular integral operator}
%%%
Let $\Gamma$ be a  rectifiable Jordan curve. We provide $\Gamma$ with the
counter-clockwise orientation. The curve $\Gamma$ divides the complex plane
$\C$ into a bounded connected component $D^+$ and an unbounded
connected component $D^-$. Without loss of generality we suppose that
$0\in D^+$. Let $X$ be a Banach function space and $w:\Gamma\to[0,\infty]$
be a weight. Then the weighted Banach function space $X_w$
is a linear normed space which becomes a Banach function space
whenever $w\in X$ and $1/w\in X'$ (see Lemma~\ref{le:WBFS}).
%%%%%%%%%%%%%%%%%%%%%%%%%%%%%%%%%%%%%%%%%%%%%%%%%%%%%%%%%%%%%%%%%%%%%%%%
\begin{theorem}\label{th:S}
Let $\Gamma$ be a  rectifiable Jordan curve, let $w:\Gamma\to[0,\infty]$
be a weight, and let $X$ be a Banach function space. If the Cauchy
singular integral operator $S$ is bounded on the weighted Banach
function space $X_w$, then $w\in A_X(\Gamma)$.
\end{theorem}
%%%%%%%%%%%%%%%%%%%%%%%%%%%%%%%%%%%%%%%%%%%%%%%%%%%%%%%%%%%%%%%%%%%%%%%%
This theorem was proved for weighted rearrangement-invariant
Banach function spaces in a slightly different form in
\cite[Theorem~3.2]{K98} (see also \cite[Theorem~4.3]{K96}
and \cite[Theorem~4.8]{BK97}).
First, as in \cite[Lemma~3.3]{K98}, by using the Landau lemma
for the Banach function space $X$ (see \cite[Ch.~1, Lemma~2.7]{BS88}),
we show that $w\in X$ and $1/w\in X'$. Then, by Lemma~\ref{le:WBFS}(b),
the weighted Banach function space $X_w$ is itself a Banach function space.
The proof of \cite[Theorem~3.2]{K98} (see also \cite[Section~3]{Dissertation})
does not use the rearrangement-invariant property of the space $X$,
so it works for arbitrary weighted Banach function spaces.

The question about the sufficiency of the condition $w\in A_X(\Gamma)$
for the boundedness of the Cauchy singular integral operator $S$
on weighted Banach function spaces $X_w$ is open. We know only that
this condition is sufficient for the boundedness in the case
of Lebesgue spaces $X=L^p, 1<p<\infty$, that is, when
$A_X(\Gamma)=A_p(\Gamma)$ is the Muckenhoupt class (see, e.g.,
\cite[Theorem~4.15]{BK97}).

However, criteria for the boundedness of $S$ on Nakano spaces
with Khvedelidze weights $L^{p(\cdot)}_\varrho$ were recently
proved by V.~M.~Kokilashvili and S.~G.~Samko \cite{KS02-3}
under the condition that the contour $\Gamma$ is sufficiently nice.
%%%%%%%%%%%%%%%%%%%%%%%%%%%%%%%%%%%%%%%%%%%%%%%%%%%%%%%%%%%%%%%%%%%%%%%%
\begin{theorem}\label{th:Samko-Kokilashvili}
{\rm (see \cite[Theorem~2]{KS02-3}).}
Let $\Gamma$ be either a Lyapunov Jordan curve or a Radon Jordan curve
without cusps, let $\varrho$ be a Khvedelidze weight
{\rm (\ref{eq:Khvedelidze-def})}, and let $p\in\cP$. The Cauchy singular
integral operator $S$ is bounded on the weighted Nakano space
$L^{p(\cdot)}_\varrho$ if and only if
%%%
\begin{equation}\label{eq:Samko-Kokilashvili}
0<\frac{1}{p(\tau_k)}+\lambda_k<1
\quad\mbox{for all}\quad k\in\{1,\dots,n\}.
\end{equation}
\end{theorem}
%%%%%%%%%%%%%%%%%%%%%%%%%%%%%%%%%%%%%%%%%%%%%%%%%%%%%%%%%%%%%%%%%%%%%%%%
For weighted Lebesgue spaces $L^p_\varrho$ this result is classic,
for Lyapunov curves it was proved by B.~V.~Khvedelidze \cite{Khvedelidze56}
and for Radon curves without cusps by I.~I.~Danilyuk and V.~Yu.~Shelepov
\cite[Theorem~2]{DS67}. The proofs and history can be found in
\cite{Danilyuk75,GK92,Khvedelidze75,MP86}.
%%%%%%%%%%%%%%%%%%%%%%%%%%%%%%%%%%%%%%%%%%%%%%%%%%%%%%%%%%%%%%%%%%%%%%%%
\subsection{Singular integral operators}\label{section-SIO}
In the following we will assume that $\Gamma$ is a rectifiable Jordan
curve, $X$ is a Banach function space, $w:\Gamma\to[0,\infty]$
is a weight such that
\begin{enumerate}
\item[{\rm (B)}] the Cauchy singular integral operator $S$
is bounded on the weighted Banach function space $X_w$;
\item[{\rm (R)}] the weighted Banach function space $X_w$ is reflexive.
\end{enumerate}

Axiom (B) guarantees that, by Theorem~\ref{th:S}, $w\in A_X(\Gamma)$.
Therefore, $w\in X$ and $1/w\in X'$. Hence, $X_w$ is a Banach
function space with the associate space $X_{1/w}'$ and
\[
L^\infty\subset X_w\subset L^1.
\]
On the other hand, if $w\in A_X(\Gamma)$, then $\Gamma$ is
a Carleson curve. Axiom
(R) implies that the Banach dual $(X_w)^*$ of $X_w$ coincides
with its associate space $X_{1/w}'$ and the set $\cR$ of
all rational functions without poles on $\Gamma$
is dense in both $X_w$ and $X_{1/w}'$ (for details, see
Subsection~\ref{subsection:WBFS}).

The above mentioned properties of weighted Banach functions spaces
satisfying axioms (B) and (R) allow us to prove the following
statements as in the case of weighted Lebsegue spaces (see, e.g.,
\cite[Ch.~1]{GK92} and \cite[Ch.~6]{BK97}). Detailed proofs can be
found in \cite[Ch.~2]{Dissertation} (see also \cite{K98,K00})
for weighted rearrangement-invariant Banach function spaces $X_w$.
Note that  the assumption that $X$ is rearrangement-invariant is not
essential and can be omitted there.

We denote by $\cK(X_w)$ the closed two-sided ideal of all compact
operators on $X_w$ in the Banach algebra $\cB(X_w)$ of all bounded
linear operators on $X_w$. As usual, $I$ is the identity operator on
$X_w$ and $aI$ denotes the operator of multiplication by a measurable
function $a:\Gamma\to\C$.
%%%%%%%%%%%%%%%%%%%%%%%%%%%%%%%%%%%%%%%%%%%%%%%%%%%%%%%%%%%%%%%%%%%%%%%%
\begin{lemma}\label{le:mult-boundedness}
If $a\in L^\infty$, then $aI\in\cB(X_w)$ and $\|aI\|_{\cB(X_w)}\le \|a\|_\infty$.
\end{lemma}
%%%%%%%%%%%%%%%%%%%%%%%%%%%%%%%%%%%%%%%%%%%%%%%%%%%%%%%%%%%%%%%%%%%%%%%%
\begin{lemma}\label{le:proj-boundedness}
The operators
\[
P_+:=(I+S)/2,\quad P_-:=(I-S)/2
\]
are bounded projections on both $X_w$ and $X_{1/w}'$.
\end{lemma}
%%%%%%%%%%%%%%%%%%%%%%%%%%%%%%%%%%%%%%%%%%%%%%%%%%%%%%%%%%%%%%%%%%%%%%%%
\begin{lemma}\label{le:compactness-commutator}
If $a\in C$, then $aS-SaI\in\cK(X_w)$.
\end{lemma}
%%%%%%%%%%%%%%%%%%%%%%%%%%%%%%%%%%%%%%%%%%%%%%%%%%%%%%%%%%%%%%%%%%%%%%%%
On the weighted Banach function space $X_w$ (or on its dual
$(X_w)^*=X_{1/w}'$) define the operator $H_\Gamma$ by
$(H_\Gamma\varphi)(\tau):=e^{-i\theta_\Gamma(\tau)}\overline{\varphi(\tau)}$.
Note that the operator $H_\Gamma$ is additive but
$H_\Gamma(\alpha\varphi)=\overline{\alpha}\cdot H_\Gamma\varphi$ for
$\alpha\in\C$. Evidently, $H_\Gamma^2=I$.
%%%%%%%%%%%%%%%%%%%%%%%%%%%%%%%%%%%%%%%%%%%%%%%%%%%%%%%%%%%%%%%%%%%%%%%%
\begin{lemma}\label{le:S-adjoint}
The adjoint of $S\in\cB(X_w)$ is $S^*=-H_\Gamma SH_\Gamma\in\cB(X_{1/w}')$.
\end{lemma}
%%%%%%%%%%%%%%%%%%%%%%%%%%%%%%%%%%%%%%%%%%%%%%%%%%%%%%%%%%%%%%%%%%%%%%%%
For $a\in L^\infty$, put
\[
T_a:=P_+aP_++P_-, 
\quad
R_a:=aP_++P_-.
\]
%%%%%%%%%%%%%%%%%%%%%%%%%%%%%%%%%%%%%%%%%%%%%%%%%%%%%%%%%%%%%%%%%%%%%%%%
\begin{lemma}\label{le:RT}
Let $a\in L^\infty$. If one of the operators $T_a,R_a$ is semi-Fredholm,
Fredholm, left-invertible, right-invertible, invertible, then the second
operator has the same property. If the operators $T_a$ and $R_a$ are
semi-Fredholm, then
\[
n(T_a)=n(R_a),\quad d(T_a)=d(R_a).
\]
\end{lemma}
%%%%%%%%%%%%%%%%%%%%%%%%%%%%%%%%%%%%%%%%%%%%%%%%%%%%%%%%%%%%%%%%%%%%%%%%
\begin{proof}
By Lemmas~\ref{le:mult-boundedness}--\ref{le:proj-boundedness},
the operators $aI$ and $P_\pm$ are bounded on $X_w$.
The rest follows from \cite[Lemma~1.21]{KS01}.
\end{proof}

So, it is sufficient to study only one of the operators $T_a, R_a$.
We will formulate our main results for the operator $R_a$. This operator
is usually called a \textit{singular integral operator with the coefficient}
$a$. It is well known that Fredholm properties of this operator are closely 
connected with the solvability of the Riemann-Hilbert boundary value problem 
(see, e.g., \cite{CG81,GK92,LS87}).
%%%%%%%%%%%%%%%%%%%%%%%%%%%%%%%%%%%%%%%%%%%%%%%%%%%%%%%%%%%%%%%%%%%%%%%%
\subsection{Hardy type subspaces}\label{section-Hardy}
In view of Lemma~\ref{le:proj-boundedness}, one can define the following
subspaces of $X_w$:
\[
(X_w)_+:=P_+X_w,
\quad
(X_w)_-^0:=P_-X_w,
\quad
(X_w)_-:=(X_w)_-^0\dot{+}\C;
\]
the corresponding subspaces $(X_{1/w}')_+, (X_{1/w}')_-^0, (X_{1/w}')_-$
of $X_{1/w}'$ are defined analogously. Also put
%%%
\begin{eqnarray*}
L^1_+ &:=& \left\{ f\in L^1:\quad\int_\Gamma f(\tau)\tau^n d\tau=0
\quad\mbox{for}\quad n\ge 0\right\},
\\
(L^1)_-^0 &:=& \left\{ f\in L^1:\quad\int_\Gamma f(\tau)\tau^n d\tau=0
\quad\mbox{for}\quad n<0\right\},
\\
L^1_- &:=& (L^1)_-^0 \dot{+}\C.
\end{eqnarray*}
%%%%%%%%%%%%%%%%%%%%%%%%%%%%%%%%%%%%%%%%%%%%%%%%%%%%%%%%%%%%%%%%%%%%%%%%
\begin{lemma}\label{le:Privalov}
{\rm (see \cite[pp.~202--206]{Privalov50})}. We have
$L_+^1\cap (L^1)_-^0=\{0\}$ and  $L^1_+\cap L^1_-=\C$.
\end{lemma}
%%%%%%%%%%%%%%%%%%%%%%%%%%%%%%%%%%%%%%%%%%%%%%%%%%%%%%%%%%%%%%%%%%%%%%%%
\begin{lemma}\label{le:subspaces-properties}
{\rm (a)} If $f\in (X_w)_\pm$ and $g\in(X_{1/w}')_\pm$, then $fg\in L^1_\pm$.
If, in addition, $f\in (X_w)_-^0$ or $g\in (X_{1/w}')_-^0$, then
$fg\in (L^1)_-^0$.

\noindent
{\rm (b)} We have,
\[
(X_w)_+=L^1_+\cap X_w,\quad
(X_w)_-^0=(L^1)_-^0\cap X_w,\quad
(X_w)_-=L^1_-\cap X_w.
\]
\end{lemma}

This lemma is proved by analogy with \cite[Corollary~6.8]{BK97}
and \cite[Lemma~6.11]{BK97}. Here we essentially use Cauchy's theorem,
H\"older's inequality for the weighted Banach function space $X_w$,
and the density of $\cR$ in $X_w$ and in $X_{1/w}'$
(see Corollary~\ref{co:R-density}).
%%%%%%%%%%%%%%%%%%%%%%%%%%%%%%%%%%%%%%%%%%%%%%%%%%%%%%%%%%%%%%%%%%%%%%%%
\begin{lemma}\label{le:Grudsky}
Suppose $f_\pm$ is analytic in $D^\pm$ and continuous on $D^\pm\cup\Gamma$
with the possible exception of finitely many points $t_1,\dots,t_m\in\Gamma$.
Suppose that $f_\pm|\Gamma\in X_w$ and that $f_\pm$ admits the estimate
\[
|f_\pm(z)|\le M|z-t_k|^{-\mu}
\quad (k=1.\dots,m)
\]
with some $M>0,\mu>0$ for all $z\in D^\pm$ sufficiently close
to $t_k$. Then $f_\pm\in (X_w)_\pm$.
\end{lemma}
%%%%%%%%%%%%%%%%%%%%%%%%%%%%%%%%%%%%%%%%%%%%%%%%%%%%%%%%%%%%%%%%%%%%%%%%
This result goes back to S.~Grudsky \cite[Proposition~1.5]{Grudsky85}
for Lebesgue spaces. To prove this statement, we should repeat the
proof of \cite[Lemma~6.10]{BK97}, replacing $L^p(\Gamma,w)$ by
$X_w$ and using Lemma~\ref{le:subspaces-properties}. For $\mu\in(0,1]$
and Lebesgue spaces this result was known for a long time
\cite[Ch.~2, Theorem~4.8]{GK92}. We remark that for our purposes
(see Lemma~\ref{le:fact-sufficiency}) we really need this analog of
Grudsky's lemma allowing also the case $\mu>1$.
%%%%%%%%%%%%%%%%%%%%%%%%%%%%%%%%%%%%%%%%%%%%%%%%%%%%%%%%%%%%%%%%%%%%%%%%
\subsection{Two basic theorems}
Let $GL^\infty$ denote the set of all functions in $L^\infty$
which are invertible in $L^\infty$, that is, the set of functions
$a\in L^\infty$ such that
\[
\operatornamewithlimits{ess\,inf}_{\tau\in\Gamma}|a(\tau)|>0.
\]
%%%%%%%%%%%%%%%%%%%%%%%%%%%%%%%%%%%%%%%%%%%%%%%%%%%%%%%%%%%%%%%%%%%%%%%%
\begin{theorem}\label{th:necessity-semi}
Let $a,b\in L^\infty$. If the operator $aP_++bP_-$ is semi-Fredholm in $X_w$,
then $a,b\in GL^\infty$.
\end{theorem}
%%%%%%%%%%%%%%%%%%%%%%%%%%%%%%%%%%%%%%%%%%%%%%%%%%%%%%%%%%%%%%%%%%%%%%%%
\begin{theorem}\label{th:Coburn-Simonenko}
If $a\in GL^\infty$, then $\min\{n(R_a),d(R_a)\}=0$.
\end{theorem}
%%%%%%%%%%%%%%%%%%%%%%%%%%%%%%%%%%%%%%%%%%%%%%%%%%%%%%%%%%%%%%%%%%%%%%%%
Theorem~\ref{th:Coburn-Simonenko} was proved by L.~Coburn
\cite{Coburn66} for Toeplitz operators on $L^2(\T)$. In the form
presented here Theorems~\ref{th:necessity-semi} and~\ref{th:Coburn-Simonenko}
were proved by I. B. Simonenko in \cite{Simonenko68} for Lebesgue spaces
with Khvedelidze weights over Lyapunov curves. For a detailed discussion
of these theorems for weighted Lebesgue spaces, see \cite[Section~6.6]{BK97}
and \cite[Sections~7.4 and 7.5]{GK92}. In our case the proofs
are developed analogously on the basis of the results of
Subsections~\ref{section-SIO}--\ref{section-Hardy}
and the Lusin-Privalov theorem (see, e.g., \cite[p.~292]{Privalov50}).
%%%%%%%%%%%%%%%%%%%%%%%%%%%%%%%%%%%%%%%%%%%%%%%%%%%%%%%%%%%%%%%%%%%%%%%%
\subsection{The local principle of Simonenko type}\label{sec:local}
Two functions $a,b\in L^\infty$ are said to be \textit{locally equivalent
at a point} $t\in\Gamma$ if
\[
\inf\Big\{\|(a-b)c\|_\infty:\quad c\in C,\quad c(t)=1\Big\}=0.
\]
%%%%%%%%%%%%%%%%%%%%%%%%%%%%%%%%%%%%%%%%%%%%%%%%%%%%%%%%%%%%%%%%%%%%%%%%
\begin{theorem}\label{th:local-principle}
Let $a\in L^\infty$. Suppose for each $t\in\Gamma$ we are given a function
$a_t\in L^\infty$ which is locally equivalent to $a$ at $t$. If the operators
$R_{a_t}$ are Fredholm in $X_w$ for all $t\in\Gamma$, then $R_a$ is Fredholm
in $X_w$.
\end{theorem}
%%%%%%%%%%%%%%%%%%%%%%%%%%%%%%%%%%%%%%%%%%%%%%%%%%%%%%%%%%%%%%%%%%%%%%%%
For weighted Lebesgue spaces, this theorem is known as Simonenko's local
principle \cite{Simonenko65}. More information about localization techniques
can be found, e.g., in \cite{BK97,BS90,GK92,Krupnik87}.
Theorem~\ref{th:local-principle} can be proved similarly to
\cite[Theorem~6.30]{BK97} with the help of
Lemmas~\ref{le:compactness-commutator} and~\ref{le:RT}.
%%%%%%%%%%%%%%%%%%%%%%%%%%%%%%%%%%%%%%%%%%%%%%%%%%%%%%%%%%%%%%%%%%%%%%%%
\subsection{Wiener-Hopf factorization}
We say that a function $a\in L^\infty$ admits a \textit{Wiener-Hopf
factorization in the weighted Banach function space} $X_w$ if
$1/a\in L^\infty$ and $a$ can be written in the form
%%%
\begin{equation}\label{eq:WH}
a(t)=a_-(t)t^\kappa a_+(t)
\quad\mbox{a.e. on}\ \Gamma,
\end{equation}
%%%
where $\kappa\in\Z$, and the factors $a_\pm$ enjoy the following properties:
%%%
\begin{enumerate}
\item[{\rm (i)}]
$a_-\in (X_w)_-, \quad
1/a_-\in (X_{1/w}')_-,\quad
a_+\in (X_{1/w}')_+,\quad
1/a_+\in(X_w)_+,
$
\item[{\rm (ii)}]
the operator $(1/a_+)Sa_+I$ is bounded on $X_w$.
\end{enumerate}
One can prove that the number $\kappa$ is uniquely determined.
%%%%%%%%%%%%%%%%%%%%%%%%%%%%%%%%%%%%%%%%%%%%%%%%%%%%%%%%%%%%%%%%%%%%%%%%
\begin{theorem}\label{th:factorization}
A function $a\in L^\infty$ admits a Wiener-Hopf factorization {\rm (\ref{eq:WH})}
in the reflexive weighted Banach function space $X_w$ if and only if the
operator $R_a$ is Fredholm in $X_w$. If $R_a$ is Fredholm, then
its index is equal to $-\kappa$.
\end{theorem}

This theorem goes back to I.~B.~Simonenko \cite{Simonenko64,Simonenko68}.
For more about this topic we refer to \cite[Section~6.12]{BK97},
\cite[Section~5.5]{BS90}, \cite[Section~8.3]{GK92} and also to \cite{CG81,LS87}
in the case of weighted Lebesgue spaces. Simonenko's result was generalized
by the author to the case of reflexive Orlicz spaces \cite[Theorem~5.6]{K96}
and to the case of reflexive rearrangement-invariant spaces
\cite[Theorem~6.10]{K98}. In the case of reflexive weighted Banach function
spaces the proof is developed by analogy. The proof is essentially based
on the density of $\cR$ in $X_w$ and in $X_{1/w}'$,
Lemmas~\ref{le:Privalov}--\ref{le:subspaces-properties}, and
Theorems~\ref{th:necessity-semi}--\ref{th:Coburn-Simonenko}.
Detailed proofs for the results of this section can be found in
\cite[Ch.~2]{Dissertation} for weighted Banach function spaces $X_w$
provided $X$ is rearrangement-invariant. Let us remind that
this assumption can be simply omitted.
%%%%%%%%%%%%%%%%%%%%%%%%%%%%%%%%%%%%%%%%%%%%%%%%%%%%%%%%%%%%%%%%%%%%%%%%
\section{Fredholmness of singular integral operators\\
in weighted Banach function spaces}\label{sec:PC}
\subsection{Local representatives}
%%%
Fix $t\in\Gamma$. For a function $a\in PC\cap GL^\infty$ we construct a
``canonical'' function $g_{t,\gamma}$ which is locally equivalent to $a$
at the point $t\in\Gamma$. The interior and the exterior of the unit circle
can be conformally mapped onto $D^+$ and $D^-$ of $\Gamma$, respectively,
so that the point $1$ is mapped to $t$, and the points $0\in D^+$ and
$\infty\in D^-$ remain fixed. Let $\Lambda_0$ and $\Lambda_\infty$
denote the images of $[0,1]$ and $[1,\infty)\cup\{\infty\}$ under this map.
The curve $\Lambda_0\cup\Lambda_\infty$ joins $0$ to $\infty$ and
meets $\Gamma$ at exactly one point, namely $t$. Let $\arg z$ be a
continuous branch of argument in $\C\setminus(\Lambda_0\cup\Lambda_\infty)$.
For $\gamma\in\C$, define the function $z^\gamma:=|z|^\gamma e^{i\gamma\arg z}$,
where $z\in\C\setminus(\Lambda_0\cup\Lambda_\infty)$. Clearly, $z^\gamma$
is an analytic function in $\C\setminus(\Lambda_0\cup\Lambda_\infty)$. The
restriction of $z^\gamma$ to $\Gamma\setminus\{t\}$ will be denoted by
$g_{t,\gamma}$. Obviously, $g_{t,\gamma}$ is continuous and nonzero on
$\Gamma\setminus\{t\}$.

Since $a(t\pm 0)\ne 0$, we can define $\gamma_t=\gamma\in\C$ by the formulas
%%%
\begin{equation}\label{eq:local-representative}
\operatorname{Re}\gamma_t:=\frac{1}{2\pi}\arg\frac{a(t-0)}{a(t+0)},
\quad
\operatorname{Im}\gamma_t:=-\frac{1}{2\pi}\log\left|\frac{a(t-0)}{a(t+0)}\right|,
\end{equation}
%%%
where we can take any value of $\arg(a(t-0)/a(t+0))$, which implies that
any two choices of $\operatorname{Re}\gamma_t$ differ by an integer only.
Clearly, there is a constant $c_t\in\C\setminus\{0\}$ such that
$a(t\pm 0)=c_tg_{t,\gamma_t}(t\pm 0)$, which means that $a$ is locally
equivalent to $c_tg_{t,\gamma_t}$ at the point $t\in\Gamma$.
%%%%%%%%%%%%%%%%%%%%%%%%%%%%%%%%%%%%%%%%%%%%%%%%%%%%%%%%%%%%%%%%%%%%%%%%%
\subsection{Sufficient conditions for factorability of the local representative}
\begin{lemma}\label{le:fact-sufficiency}
If, for some $k\in\Z$ and $\gamma\in\C$, the operator
$\varphi_{t,k-\gamma}S\varphi_{t,\gamma-k}I$ is bounded on
the weighted Banach function space $X_w$, then
%%%
\begin{equation}\label{eq:fact-sufficiency-1}
g_{t,\gamma}(\tau)=(1-t/\tau)^{k-\gamma}\tau^k(\tau-t)^{\gamma-k},
\quad \tau\in\Gamma\setminus\{t\}
\end{equation}
%%%
is a Wiener-Hopf factorization of the function $g_{t,\gamma}$ in $X_w$.
\end{lemma}
%%%%%%%%%%%%%%%%%%%%%%%%%%%%%%%%%%%%%%%%%%%%%%%%%%%%%%%%%%%%%%%%%%%%%%%%%%%
\begin{proof}
Since the operator $\varphi_{t,k-\gamma}S\varphi_{t,k-\gamma}^{-1}I$
is bounded on $X_w$, the operator $S$ is bounded on the weighted
Banach function space $X_{\varphi_{t,k-\gamma}w}$.
By Theorem~\ref{th:S}, $\varphi_{t,k-\gamma}w\in A_X(\Gamma)$.
In that case $\Gamma$ is a Carleson curve and
$\varphi_{t,k-\gamma}w\in X$, whence $\varphi_{t,k-\gamma}\in X_w$.

Let us show that $(\tau-t)^{k-\gamma}\in (X_w)_+$. The function
$f(z):=(z-t)^{k-\gamma}$ is analytic in $D^+$ and continuous
on $D^+\cup(\Gamma\setminus\{t\})$. For $z\in D^+$,
\[
|f(z)|=|(z-t)^{k-\gamma}|
=
|z-t|^{k-\operatorname{Re}\gamma-\Theta_t(z)\operatorname{Im}\gamma},
\]
where $\Theta_\tau(z):=\arg(z-t)/(-\log|z-t|)$. As in \cite[Theorem~7.7]{K98}
and \cite[Lemma~7.1]{BK97} with the help of Lemma~\ref{le:Seifullaev} one can
show that there is a constant $M_t\in(0,\infty)$ such that
\[
|f(z)|\le |z-t|^{k-\operatorname{Re}\gamma}e^{M_t|\operatorname{Im}\gamma|(-\log|z-t|)}
=
|z-t|^{k-\operatorname{Re}\gamma-M_t|\operatorname{Im}\gamma|}
\]
for all $z$ in a small neighborhood of $t$. By Lemma~\ref{le:Grudsky},
$(\tau-t)^{k-\gamma}\in (X_w)_+$. Analogously one can prove that
\[
(\tau-t)^{\gamma-k}\in (X_{1/w}')_+,
\quad
(1-t/\tau)^{k-\gamma}\in (X_w)_-,
\quad
(1-t/\tau)^{\gamma-k}\in (X_{1/w}')_-.
\]
These facts together with
the boundedness of $\varphi_{t,k-\gamma}S\varphi_{t,\gamma-k}I$
on the space $X_w$ show that (\ref{eq:fact-sufficiency-1})
is indeed a Wiener-Hopf factorization of the function $g_{t,\gamma}$.
\end{proof}
%%%%%%%%%%%%%%%%%%%%%%%%%%%%%%%%%%%%%%%%%%%%%%%%%%%%%%%%%%%%%%%%%%%%%%%%%%%
\subsection{Necessary conditions for factorability of the local representative}
\begin{theorem}\label{th:fact-necessity}
If the function $g_{t,\gamma}$ admits a Wiener-Hopf factorization in
the weighted Banach function space $X_w$, then
%%%
\begin{equation}\label{eq:fact-necessity-1}
-\operatorname{Re}\gamma+\theta\alpha_t^*(-\operatorname{Im}\gamma)
+(1-\theta)\beta_t^*(-\operatorname{Im}\gamma)\notin\Z
\end{equation}
%%%
for all $\theta\in[0,1]$. Moreover, there exists an $l\in\Z$ such that
$\varphi_{t,l-\gamma}w\in A_X(\Gamma)$.
\end{theorem}
%%%%%%%%%%%%%%%%%%%%%%%%%%%%%%%%%%%%%%%%%%%%%%%%%%%%%%%%%%%%%%%%%%%%%%%%%%%
\begin{proof}
The idea of the proof (in the case of weighted Lebesgue spaces)
goes back to I.~Spitkovsky \cite{Spitkovsky92} and it was further developed
by A.~B\"ottcher and Yu.~I.~Karlovich \cite[Proposition~7.2]{BK97}.
This idea was applied to the proof in the case of reflexive rearrangement-invariant
Banach function spaces (with weights) in \cite[Theorem~7.6]{K98}
and \cite[Theorem~4.1]{K00}. Since, for our (more general) case, the arguments
are the same, we point out only the main steps.

By Theorem~\ref{th:factorization}, the operator $g_{t,\gamma}P_++P_-$
is Fredholm. Then there exists a $c>0$ such that the operators
$g_{t,\gamma-\eps}P_++P_-$ are Fredholm for all $\eps\in(-c,c)$.
Applying Theorem~\ref{th:factorization} again, we infer that all functions
$g_{t,\gamma-\eps}$ admit a Wiener-Hopf factorization in $X_w$. By using its
definition, one can show that there exists an $l\in\Z$ such that the
operators $\varphi_{t,l-\gamma+\eps}S\varphi_{t,l-\gamma+\eps}^{-1}I$
are bounded on $X_w$ for all $\eps\in(-c,c)$. In that case, by
Theorem~\ref{th:S}, $\varphi_{t,l-\gamma+\eps}w\in A_X(\Gamma)\subset A_X(\Gamma,t)$.
By Lemma~\ref{le:Q-reg},
%%%
\begin{equation}\label{eq:fact-necessity-2}
0\le(Q_t(\varphi_{t,l-\gamma+\eps}w))\le\beta(Q_t(\varphi_{t,l-\gamma+\eps}w))\le 1.
\end{equation}
%%%
From Lemma~\ref{le:indicator*} and (\ref{eq:fact-necessity-2}) it follows
that
\[
0\le l+\eps-\operatorname{Re}\gamma+\alpha_t^*(-\operatorname{Im}\gamma)
\le
 l+\eps-\operatorname{Re}\gamma+\beta_t^*(-\operatorname{Im}\gamma)\le 1
\]
for all $\eps\in(-c,c)$. Hence,
\[
-l<-\operatorname{Re}\gamma+\theta\alpha_t^*(-\operatorname{Im}\gamma)
+(1-\theta)\beta_t^*(-\operatorname{Im}\gamma)<l-1
\]
for every $\theta\in[0,1]$. Thus, (\ref{eq:fact-necessity-1}) holds for every
$\theta\in[0,1]$.
\end{proof}
%%%%%%%%%%%%%%%%%%%%%%%%%%%%%%%%%%%%%%%%%%%%%%%%%%%%%%%%%%%%%%%%%%%%%%%%%%%
\subsection{Necessary conditions for Fredholmness}
%%%
Now we are in a position to state the main result of this paper.
%%%%%%%%%%%%%%%%%%%%%%%%%%%%%%%%%%%%%%%%%%%%%%%%%%%%%%%%%%%%%%%%%%%%%%%%%%
\begin{theorem}\label{th:necessity}
Let $\Gamma$ be a rectifiable Jordan curve, let $w:\Gamma\to[0,\infty]$
be a weight, and let $X$ be a Banach function space. Suppose the Cauchy
singular integral operator $S$ is bounded on the weighted
Banach function space $X_w$ and $X_w$ is reflexive.
If the operator $aP_++P_-$, where $a\in PC$, is Fredholm in $X_w$, then
$a\in GL^\infty$ and
\begin{eqnarray}
\label{eq:necessity-1}
&-&\frac{1}{2\pi}\arg\frac{a(t-0)}{a(t+0)}
\\
&+&
\theta\alpha_t^*\left(\frac{1}{2\pi}
\log\left|\frac{a(t-0)}{a(t+0)}\right|\right)
+
(1-\theta)\beta_t^*\left(
\frac{1}{2\pi}\log\left|\frac{a(t-0)}{a(t+0)}\right|\right)
\not\in\Z
\nonumber
\end{eqnarray}
for all $t\in\Gamma$ and all $\theta\in[0,1]$.
\end{theorem}
%%%%%%%%%%%%%%%%%%%%%%%%%%%%%%%%%%%%%%%%%%%%%%%%%%%%%%%%%%%%%%%%%%%%%%%%%
\begin{proof}
The proof is developed by analogy with the proof of necessity part of
\cite[Theorem~7.8]{K98} (see also \cite[Proposition~7.3]{BK97}).

If $R_a$ is Fredholm, then, by Theorem~\ref{th:necessity-semi},
$a\in GL^\infty$. Fix an arbitrary $t\in\Gamma$. Choose $\gamma=\gamma_t\in\C$
as in (\ref{eq:local-representative}).
Then the function $a$ is locally equivalent to $c_tg_{t,\gamma_t}$ at
the point $t$, where $c_t\in\C\setminus\{0\}$ is some constant. If
$\tau\in\Gamma\setminus\{t\}$, then $g_{t,\gamma_t}$ is continuous
and nonzero at $\tau$. Hence, it is locally equivalent to the
nonzero constant $b_\tau:=g_{t,\gamma_t}(\tau)$ at $\tau$.
Clearly, the operator $R_{b_\tau}:=b_\tau P_++P_-$ is invertible,
$(b_\tau P_++P_-)^{-1}=b_\tau^{-1}P_++P_-$. Therefore, the operator
$R_{b_\tau}$ is Fredholm for every $\tau\in\Gamma\setminus\{t\}$.
Remind that the function $g_{t,\gamma_t}$ is locally equivalent
to the function $c_t^{-1}a$. Since
%%%
\begin{equation}\label{eq:necessity-3}
R_{c_t^{-1}}R_a=P_+c_t^{-1}aP_++P_-=T_{c_t^{-1}a}
\end{equation}
%%%
and the operator $R_{c_t^{-1}}$ is invertible, from Lemma~\ref{le:RT}
and (\ref{eq:necessity-3}) it follows that $R_a$ is Fredholm if and only
if $R_{c_t^{-1}a}$ is Fredholm. Therefore, applying
Theorem~\ref{th:local-principle}, we infer that the operator $R_{g_{t,\gamma_t}}$
is Fredholm. By Theorem~\ref{th:factorization}, the function $g_{t,\gamma_t}$
admits a Wiener-Hopf factorization in $X_w$. From Theorem~\ref{th:fact-necessity}
it follows that
%%%
\begin{equation}\label{eq:necessity-4}
-\operatorname{Re}\gamma_t
+\theta\alpha_t^*(-\operatorname{Im}\gamma_t)
+(1-\theta)\beta_t^*(-\operatorname{Im}\gamma_t)
\notin\Z
\end{equation}
%%%
for all $\theta\in[0,1]$. Since $t\in\Gamma$ is arbitrary, from
(\ref{eq:local-representative}) and (\ref{eq:necessity-4}) we conclude
that (\ref{eq:necessity-1}) holds for every $t\in\Gamma$ and every
$\theta\in[0,1]$.
\end{proof}
%%%%%%%%%%%%%%%%%%%%%%%%%%%%%%%%%%%%%%%%%%%%%%%%%%%%%%%%%%%%%%%%%%%%%%%%%
\subsection{Lower estimates for essential norms}
For an operator $A\in\cB(X_w)$, let
\[
|A|_{X_w}:=\inf_{K\in\cK(X_w)}\|A+K\|_{\cB(X_w)}
\]
be its \textit{essential norm} in $X_w$.
%%%%%%%%%%%%%%%%%%%%%%%%%%%%%%%%%%%%%%%%%%%%%%%%%%%%%%%%%%%%%%%%%%%%%%%%%
\begin{theorem}\label{th:estimate}
Let $\Gamma$ be a rectifiable Jordan curve, let $w:\Gamma\to[0,\infty]$
be a weight, and let $X$ be a Banach function space.
If the Cauchy singular integral operator $S$ is bounded on
the weighted Banach function space $X_w$ and $X_w$ is reflexive, then
\[
|S|_{X_w}\ge\cot\Big(\pi\Lambda_{\Gamma,X,w}/2 \Big),
\quad
|P_\pm|_{X_w}\ge 1/\sin(\pi\Lambda_{\Gamma,X,w}),
\]
where
\[
\Lambda_{\Gamma,X,w}:=\inf_{t\in\Gamma}
\min\Big\{\alpha(Q_tw),1-\beta(Q_tw)\Big\}.
\]
\end{theorem}

This statement is proved by a literal repetition of the proof of
\cite[Theorem~4.5]{K00} using the scheme of \cite[Ch.~9, Theorem~9.1]{GK92}. 
One can find more information about estimates
of (essential) norms on weighted Lebesgue spaces in
\cite[Ch.~13]{GK92} and \cite[Ch.~2]{Krupnik87}.
%%%%%%%%%%%%%%%%%%%%%%%%%%%%%%%%%%%%%%%%%%%%%%%%%%%%%%%%%%%%%%%%%%%%%%%%%
\section{Fredholmness of singular integral operators\\
in weighted Nakano spaces}\label{sec:PC-Nakano}
\subsection{Necessary conditions for Fredholmness}
The necessary conditions for the Fredholmness of $R_a$ in weighted Nakano
spaces have a simpler form than in the general case because we can replace
the indicator functions $\alpha_t^*$ and $\beta_t^*$
by the indicator functions $1/p(t)+\alpha_t$ and $1/p(t)+\beta_t$,
respectively. More precisely, the next theorem is true.
%%%%%%%%%%%%%%%%%%%%%%%%%%%%%%%%%%%%%%%%%%%%%%%%%%%%%%%%%%%%%%%%%%%%%%%%
\begin{theorem}\label{th:Nakano-necessity}
Let $\Gamma$ be a rectifiable Jordan curve, let $w:\Gamma\to[0,\infty]$ be a weight,
and let $p\in\cP$. Suppose the Cauchy singular integral operator is bounded on
the weighted Nakano space $L^{p(\cdot)}_w$.  If the operator $aP_++P_-$, where
$a\in PC$, is Fredholm in $L^{p(\cdot)}_w$, then $a\in GL^\infty$
and
%%%
\begin{eqnarray}
\label{eq:Nakano-necessity-1}
&-&
\frac{1}{2\pi}\arg\frac{a(t-0)}{a(t+0)}+\frac{1}{p(t)}
\\
&+&
\theta\alpha_t\left(\frac{1}{2\pi}\log\left|\frac{a(t-0)}{a(t+0)}\right|\right)
+
(1-\theta)
\beta_t\left(\frac{1}{2\pi}\log\left|\frac{a(t-0)}{a(t+0)}\right|\right)
\not\in\Z
\nonumber
\end{eqnarray}
%%%
for all $t\in\Gamma$ and all $\theta\in[0,1]$.
\end{theorem}
%%%%%%%%%%%%%%%%%%%%%%%%%%%%%%%%%%%%%%%%%%%%%%%%%%%%%%%%%%%%%%%%%%%%%%%%%
\begin{proof}
Since $p\in\cP$, inequalities (\ref{eq:p-cont-reflexive}) are satisfied.
In that case, by Lemma~\ref{le:Nakano-duality}, the non-weighted Nakano spaces
$L^{p(\cdot)}$ is reflexive. On the other hand, by Theorem~\ref{th:S},
$w\in L^{p(\cdot)}$ and $1/w\in (L^{p(\cdot)})'$. Then the weighted
Nakano space $L^{p(\cdot)}_w$ is also reflexive, due to
Corollary~\ref{co:weighted-reflexivity}. Thus, all assumptions of
Theorem~\ref{th:necessity} are satisfied and we can repeat its proof.
In view of Theorem~\ref{th:fact-necessity}, there exists an $l\in\Z$ such that
$\varphi_{t,l-\gamma_t}w\in A_{L^{p(\cdot)}}(\Gamma)$, where $\gamma_t$
is given by (\ref{eq:local-representative}). In that case, by
Lemma~\ref{le:Nakano-ind-Carleson},
%%%
\begin{eqnarray*}
&&
-\operatorname{Re}\gamma_t
+\theta\alpha_t^*(-\operatorname{Im}\gamma_t)
+(1-\theta)\beta_t^*(-\operatorname{Im}\gamma_t)
\\
&&
=
-\operatorname{Re}\gamma_t+1/p(t)
+\theta\alpha_t(-\operatorname{Im}\gamma_t)
+(1-\theta)\beta_t(-\operatorname{Im}\gamma_t).
\end{eqnarray*}
%%%
Therefore, we can replace condition (\ref{eq:necessity-1})
by condition (\ref{eq:Nakano-necessity-1}) in the case of
weighted Nakano spaces.
\end{proof}

For Lebesgue spaces $L^p_w$ with Muckenhoupt weights $w$
(that is, in the case when $p(\cdot)$ is constant),
condition (\ref{eq:Nakano-necessity-1})
becomes also sufficient for the Fredholmness of $R_a$
(see \cite[Proposition~7.3]{BK97}).
%%%%%%%%%%%%%%%%%%%%%%%%%%%%%%%%%%%%%%%%%%%%%%%%%%%%%%%%%%%%%%%%%%%%%%%%%
\subsection{Lower estimates for essential norms}
\begin{theorem}
Let $\Gamma$ be a rectifiable Jordan curve, let $w:\Gamma\to[0,\infty]$ be a weight,
and let $p\in\cP$. If the Cauchy singular integral operator is bounded on
the weighted Nakano space $L^{p(\cdot)}_w$, then
\[
|S|_{L^{p(\cdot)}_w}\ge \cot\Big(\pi\Lambda_{\Gamma,p,w}/2\Big),
\quad
|P_\pm|_{L^{p(\cdot)}_w}\ge 1/\sin(\pi\Lambda_{\Gamma,p,w}),
\]
where
\[
\Lambda_{\Gamma,p,w}:=\inf_{t\in\Gamma}
\min\left\{\frac{1}{p(t)}+\mu_t,1-\frac{1}{p(t)}-\nu_t\right\}.
\]
\end{theorem}
%%%%%%%%%%%%%%%%%%%%%%%%%%%%%%%%%%%%%%%%%%%%%%%%%%%%%%%%%%%%%%%%%%%%%%%%
By Theorem~\ref{th:S}, $w\in A_{L^{p(\cdot)}}(\Gamma)$. Therefore,
the latter theorem immediately follows from Theorem~\ref{th:estimate}
and Theorem~\ref{th:disintegration}.

If $\log w\in VMO(\Gamma)$ (in particular, if $w=1$), then from Lemma~\ref{le:VMO-indices}
and (\ref{eq:p-cont-reflexive}) it follows that
%%%
\begin{eqnarray*}
\Lambda_{\Gamma,p,w}
&=&
\inf_{t\in\Gamma}\min\left\{\frac{1}{p(t)},1-\frac{1}{p(t)}\right\}
=
\min\left\{
\inf_{t\in\Gamma}\frac{1}{p(t)},1-\sup_{t\in\Gamma}\frac{1}{p(t)}
\right\}
\\
&=&
\min\Big\{1/p^*,1-1/p_*\Big\}.
\end{eqnarray*}
%%%%%%%%%%%%%%%%%%%%%%%%%%%%%%%%%%%%%%%%%%%%%%%%%%%%%%%%%%%%%%%%%%%%%%%%%
\subsection{Fredholm criterion}
\begin{theorem}\label{th:criterion}
Let $\Gamma$ be either a Lyapunov Jordan curve or a Radon Jordan curve
without cusps, let $p\in\cP$, and let $\varrho$ be a Khvedelidze
weight {\rm (\ref{eq:Khvedelidze-def})} satisfying
{\rm (\ref{eq:Samko-Kokilashvili})}.
Then the operator $aP_++P_-$, where $a\in PC$, is Fredholm in
the weighted Nakano space $L^{p(\cdot)}_\varrho$ if and only if
%%%
\begin{equation}\label{eq:criterion-3}
a(t\pm 0)\ne 0,
\quad
-\frac{1}{2\pi}\arg\frac{a(t-0)}{a(t+0)}+ \frac{1}{p(t)}+\lambda(t)\notin\Z
\end{equation}
%%%
for all $t\in\Gamma$, where
%%%
\begin{equation}\label{eq:criterion-4}
\lambda(t)=\left\{
\begin{array}{lcl}
\lambda_k, &\mbox{if} & t=\tau_k, \quad k\in\{1,\dots,n\},\\
0,         &\mbox{if} & t\notin\Gamma\setminus\{\tau_1,\dots,\tau_n\}.
\end{array}
\right.
\end{equation}
\end{theorem}
%%%%%%%%%%%%%%%%%%%%%%%%%%%%%%%%%%%%%%%%%%%%%%%%%%%%%%%%%%%%%%%%%%%%%%%%%%
\begin{proof}
By Theorem~\ref{th:Samko-Kokilashvili}, the operator $S$ is bounded
on the (reflexive) weighted Nakano space $L^{p(\cdot)}_\varrho$.

\textit{Necessity}.
By Proposition~\ref{pr:nice-curves}, for Lyapunov curves and Radon curves without
cusps, we have $\delta_t^-=\delta_t^+=0$ whenever $t\in\Gamma$.
By Lemma~\ref{le:Nakano-ind-smooth}, the indicator
functions of the pair $(\Gamma,\varrho)$ are constants
$\alpha_t(x)=\mu_t, \beta_t(x)=\nu_t$ for $x\in\R$, where the indices of
powerlikeness $\mu_t,\nu_t$ of the Khvedelidze weight
(\ref{eq:Khvedelidze-def}) coincide with $\lambda(t)$ given by
(\ref{eq:criterion-4}). Thus,
\[
\theta\alpha_t\left(
\frac{1}{2\pi}\log\left|\frac{a(t-0)}{a(t+0)}\right|
\right)
+
(1-\theta)\beta_t\left(
\frac{1}{2\pi}\log\left|\frac{a(t-0)}{a(t+0)}\right|
\right)
=\lambda(t)
\]
for every $\theta\in[0,1]$ and every $t\in\Gamma$. Therefore, the necessity
of conditions (\ref{eq:criterion-3}) follows from Theorem~\ref{th:Nakano-necessity}.
The necessity part is proved.

\textit{Sufficiency}. From (\ref{eq:criterion-3}) it follows that
for every $t\in\Gamma$, there exists an $m_t\in\Z$ such that
\[
0<m_t-\operatorname{Re}\gamma_t+\frac{1}{p(t)}+\lambda(t)<1,
\]
where $\gamma_t$ is given by (\ref{eq:local-representative}).
By Theorem~\ref{th:Samko-Kokilashvili}, the operator $S$ is bounded on the
weighted Nakano space $L^{p(\cdot)}_{\widetilde{\varrho}_t}$, where
\[
\widetilde{\varrho}_t(\tau):=|\tau-t|^{m_t-\operatorname{Re}\gamma_t}\varrho(\tau),
\quad\tau\in\Gamma.
\]
In view of (\ref{eq:canonical-weight}) and Proposition~\ref{pr:nice-curves},
there exist constants $C_1(t), C_2(t)\in (0,\infty)$ such that
\[
C_1(t)\widetilde{\varrho}_t(\tau)
\le
\varphi_{t,m_t-\gamma_t}(\tau)
\le
C_2(t)\widetilde{\varrho}_t(\tau),
\quad\tau\in\Gamma\setminus\{t\}.
\]
Therefore, $S\in\cB(L^{p(\cdot)}_{\widetilde{\varrho}_t})$
if and only if
$\varphi_{t,m_t-\gamma_t}S\varphi_{t,\gamma_t-m_t}I\in\cB(L^{p(\cdot)}_\varrho)$.
By Lemma~\ref{le:fact-sufficiency}, the function $g_{t,\gamma_t}$ admits a
Wiener-Hopf factorization in the weighted Nakano space $L^{p(\cdot)}_\varrho$.
Due to Theorem~\ref{th:factorization}, for every $t\in\Gamma$,
the operator $g_{t,\gamma_t}P_++P_-$ is Fredholm. Then the operator
$cg_{t,\gamma_t}P_++P_-$ is Fredholm for $c\in\C\setminus\{0\}$
(see the proof of Theorem~\ref{th:necessity}).

Since the function
$c_tg_{t,\gamma_t}$ with a specially chosen constant $c_t\in\C\setminus\{0\}$
is locally equivalent to the function $a\in PC$
at every point $t\in\Gamma$, in view of Theorem~\ref{th:local-principle},
the operator $R_a=aP_++P_-$ is Fredholm in the weighted Nakano space
$L^{p(\cdot)}_\varrho$.
\end{proof}
%%%%%%%%%%%%%%%%%%%%%%%%%%%%%%%%%%%%%%%%%%%%%%%%%%%%%%%%%%%%%%%%%%%%%%%%%%
In Theorem~\ref{th:criterion} the coefficient $a$ can have a countable set
of jumps. If $a$ has only a finite number of jumps and $\varrho=1$, this
result was obtained in \cite[Theorem~A]{KS03-4} (as well as a formula for
the index of the operator $R_a$). Note that the transition from finitely
many to infinitely many jumps is more or less standard 
(see \cite[Section~9.8]{GK92} for Lebesgue spaces with Khvedelidze weights
over Lyapunov curves), using the stability
of Fredholm operators and localization techniques (see Section~\ref{sec:local}).
We give the proof of Theorem~\ref{th:criterion} here for completeness.
For Lebesgue spaces with Khvedelidze
weights over Lyapunov curves the corresponding result was obtained
in the late sixties by I.~Gohberg and N.~Krupnik \cite[Ch.~9]{GK92}.
%%%%%%%%%%%%%%%%%%%%%%%%%%%%%%%%%%%%%%%%%%%%%%%%%%%%%%%%%%%%%%%%%%%%%%%%%
\subsection*{Acknowledgments}
I would like to express my deep gratitude to Professor Lech Maligranda
(Lule\r{a} University of Technology, Sweden) for historical remarks
concerning Nakano spaces and to Professor Albrecht B\"ottcher
(Chemnitz Technical University, Germany) for useful discussions.
%%%%%%%%%%%%%%%%%%%%%%%%%%%%%%%%%%%%%%%%%%%%%%%%%%%%%%%%%%%%%%%%%%%%%%%%%

\end{document}